\documentclass[11pt,onecolumn,a4paper,english]{article}
\usepackage[pdftex]{graphicx}
\usepackage{textpos}
\usepackage{babel,varioref}
\usepackage{fancyhdr}
\usepackage{amsmath}
\usepackage{amsfonts}
\usepackage{amssymb}
\usepackage{array}
\usepackage{color}
\usepackage{colortbl}
\usepackage{amsthm}
\usepackage{subfig}
 \DeclareGraphicsExtensions{.pdf,.jpeg,.png}

\begin{document}
%
%
\title{The Generalization of the Decomposition of Functions by Energy operators (Part II) and some Applications}
%
%
\maketitle
\author{J.-P.~Montillet, $j.p.montillet@anu.edu.au$}
\begin{abstract}
\boldmath{This work introduces the families of generalized energy operators $\big ([[.]^p]_k^+\big)_{k\in\mathbb{Z}}$ and $\big ([[.]^p]_k^- \big)_{k\in\mathbb{Z}}$ ($p$ in $\mathbb{Z}^+$). One shows that with $\bold{Lemma}$ $1$, the successive derivatives of $\big ([[$f$]^{p-1}]_1^+ \big)^n$ ($n$ in $\mathbb{Z}$, $n\neq 0$) can be decomposed with the generalized energy operators  $\big ([[.]^p]_k^+\big)_{k\in\mathbb{Z}}$ when $f$ is in the subspace $\mathbf{S}_p^-(\mathbb{R})$. With $\bold{Theorem}$ $1$ and $f$ in $\mathbf{s}_p^-(\mathbb{R})$, one can decompose uniquely the successive derivatives of $\big ([[$f$]^{p-1}]_1^+ \big)^n$ ($n$ in $\mathbb{Z}$, $n\neq 0$) with the generalized energy operators  $\big ([[.]^p]_k^+\big)_{k\in\mathbb{Z}}$ and $\big ([[.]^p]_k^-\big)_{k\in\mathbb{Z}}$. $\mathbf{S}_p^-(\mathbb{R})$ and $\mathbf{s}_p^-(\mathbb{R})$ ($p$ in $\mathbb{Z}^+$) are subspaces of the Schwartz space $\mathbf{S}^-(\mathbb{R})$. These results generalize the work of \cite{JPMontillet2013}. The second fold of this work is the application of the generalized energy operator families onto the solutions of linear partial differential equations. The solutions are functions of two variables and defined in subspaces of $\mathbf{S}^-(\mathbb{R}^2)$.  The theory is then applied to the Helmholtz equation. In this specific case, the use of generalized energy operators in the general solution of this PDE extends the results of \cite{JPMontillet2010}. This work ends with some numerical examples. We also underline that this theory could possibly open some applications in astrophysics and aeronautics. 
%
%
}
\end{abstract}


%
%
%
%
%
%
%
%
\section{Introduction and guidelines}\label{Introduction part}
%
%
%
Two decades ago, an energy operator (${\Psi}_{R}^{-}$) was first defined in \cite{Kaiser90}. Since then, this work has been extensively used in signal processing (e.g., \cite{Bovik93} or \cite{Hamila1999}) and image processing \cite{Felsberg}. The bilinearity properties of this operator were studied in \cite{Boudraa et al.2009}. More recently, the author in \cite{JPMontillet2010} introduced the energy operators (${\Psi}_{R}^{+}$, ${\Psi}_{R}^{-}$) and developed a general method for separating smooth real-valued finite energy functions ($f$) in time and space with application to the wave equation. In \cite{JPMontillet2013}, the author introduced the family of differential energy operators (DEOs) $({\Psi}_{k}^{-})_{k\in\mathbb{Z}^+}$ and $({\Psi}_{k}^{+})_{k\in\mathbb{Z}^+}$, and showed the decomposition of the $s$-th  derivative of $f^n$ ($n$ in $\mathbb{Z}$, $n \neq 0$) with the DEO families. In addition, this work demonstrated some properties of the family of energy operators and the application to  the energy function $\mathcal{E}(f^n)$.
\\ This work starts with recalling the notations used with energy operators, important definitions (e.g., $\bold{Definition}$ $1$ and $\bold{Definition}$ $2$) and results ( e.g., $\bold{Lemma}$ $0$ and $\bold{Theorem}$ $0$) from the recently published work of  \cite{JPMontillet2013}. It emphasizes the notion of decomposing finite energy functions $f$ in a Schwartz space ($\mathbf{S}^-(\mathbb{R})$) with families of energy operator. In the second part, $\bold{Definition}$ $3$ defines an energy space ($\mathbf{E}$) subset of $\mathbf{S}^-(\mathbb{R})$ which is used to define other subspaces (e.g., $\mathbf{S}_p^-(\mathbb{R})$, $\mathbf{s}_p^-(\mathbb{R})$). In Section \ref{sectionenergyopdefandgen}, one allows to define families of generalized energy operators. Their properties are also shown in $\bold{Proposition}$ $3$ (e.g., bilinearity and derivative chain rule property). To extend the work in \cite{JPMontillet2013}, the statements of $\bold{Lemma}$ $0$ and $\bold{Theorem}$ $0$ are generalized with $\bold{Lemma}$ $1$ and $\bold{Theorem}$ $1$. Thus in $\bold{Lemma}$ $1$, the work emphasizes in particular $\mathbf{S}_p^-(\mathbb{R})$ a subset of $\mathbf{S}^-(\mathbb{R})$ where the decomposition of $\big ([[f]^{p-1}]_1^+ \big)^n$ ($n$ in $\mathbb{Z}$, $n\neq 0$) with the generalized energy operators  $\big ([[.]^p]_k^+\big)_{k\in\mathbb{Z}}$ is valid. Whereas in $\bold{Theorem}$ $1$, one shows the unique decomposition of $\big ([[f]^{p-1}]_1^+ \big)^n$ ($n$ in $\mathbb{Z}$, $n\neq 0$) with the generalized energy operators  $\big ([[.]^p]_k^+\big)_{k\in\mathbb{Z}}$ and $\big ([[.]^p]_k^-\big)_{k\in\mathbb{Z}}$ when $f$ is in $\mathbf{s}_p^-(\mathbb{R})$ a subspace of $\mathbf{S}^-(\mathbb{R})$. Note that in Section $2$ to $4$,  $f$  is a function of time ($\partial_t$).
\newline The second part of the work focuses on the application of the theory developed in \cite{JPMontillet2013} and in this paper with the generalized energy operators  to the linear PDEs. In this application, the solutions are functions of two variables defined onto subspaces of  $\mathbf{S}^{-}(\mathbb{R}^2)$. The definition of these subspaces include the energy operators and the generalized energy operators using $\mathbf{Theorem}$ $0$ and $\mathbf{Theorem}$ $1$. With the help of these subspaces, we define a  new model for the solutions of the linear PDEs on to $\mathbf{S}^-(\mathbb{R}^2)$. This work ends with a numerical application of this model in the case of the Helmholtz equation and the particular case of evanescent waves. We also estimate the averaged power using generalized energy operators. 
%
%
%
%
\section{Preliminaries}\label{preliminariesSection}
Throughout this work, $f^n$ for any $n$ in $\mathbb{Z}^+-\{0\}$ is supposed to be a smooth real-valued and finite energy function, and in the Schwartz space $\mathbf{S}^{-}(\mathbb{R})$ defined as:
\begin{equation}\label{SRRRRR}
\mathbf{S}^{-}(\mathbb{R}) =\{f \in \mathbf{C}^{\infty}(\mathbb{R}), \qquad {sup}_{t<0} |t^k||\partial_t^j f(t)|<\infty,\qquad \forall k \in \mathbb{Z}^+, \qquad \forall j \in \mathbb{Z}^+ \}
\end{equation}
Sometime $f^n$ can also be analytic if its development in Taylor-Series is relevant to this work. The choice of $f^n$ (for any $n$ in $\mathbb{Z}^+-\{0\}$) in the Schwartz space $\mathbf{S}^{-}(\mathbb{R})$ is based on the work developed in \cite{JPMontillet2013},  as we are dealing with multiple integrations or derivatives of $f^n$ when applying the energy operators $({\Psi}_{k}^{-})_{k\in\mathbb{Z}^+}$, $({\Psi}_{k}^{+})_{k\in\mathbb{Z}^+}$ and later on the generalized energy operators. 
\\ In the following, let us call the set $\mathcal{F}(\mathbf{S}^{-}(\mathbb{R}),\mathbf{S}^{-}(\mathbb{R}))$ all functions/operators defined such as $\gamma:$ $\mathbf{S}^{-}(\mathbb{R})$ $\rightarrow$ $\mathbf{S}^{-}(\mathbb{R})$. Let us recall some definitions and important results given in \cite{JPMontillet2013}.
\vspace{1.0em}
\newline $\bold{Definition}$ $1$: for all $f$ in $\mathbf{S}^{-}(\mathbb{R})$, for all $v\in\mathbb{Z}^+-\{0\}$, for all  $n\in\mathbb{Z}^+$ and $n>1$, the family of operators $(\Psi_k)_{k \in \mathbb{Z}}$ (with $(\Psi_k)_{k \in \mathbb{Z}}$ $\subseteq$ $\mathcal{F}(\mathbf{S}^{-}(\mathbb{R}),\mathbf{S}^{-}(\mathbb{R}))$) decomposes $\partial_t^v$$f^n$ in $\mathbb{R}$, if it exists $(N_j)_{j\in \mathbb{Z}^+ \cup \{0\}}$ $\subseteq$ $\mathbb{Z^+}$,  $(C_i)_{i=-N_j}^{N_j}$ $\subseteq$ $\mathbb{R}$, and it exists $(\alpha_j)$ and $l$ in $\mathbb{Z^+}\cup\{0\}$ (with $l<v$) 
\\ such as $\partial_t^v$$f^n = \sum_{j=0}^{v-1} \big(_{j}^{v-1} \big) \partial_t^{v-1-j} f^{n-l} \sum_{k=-N_j}^{N_j} C_k \Psi_k(\partial_t^{\alpha_k}f)$.
\vspace{1.0em}
\\ Note that Definition $1$ varies slightly from \cite{JPMontillet2013} (see updated version in Arxiv \cite{JPMontillet2013}). In addition, one has to define $\mathbf{s}^{-}(\mathbb{R})$ as:
\begin{equation}
\mathbf{s}^{-}(\mathbb{R}) = \{ f \in  \mathbf{S}^{-}(\mathbb{R}) | f \notin  (\cup_{k \in \mathbb{Z}} Ker(\Psi^{+}_k))\cup(\cup_{k \in \mathbb{Z}-\{1\}} Ker(\Psi^{-}_k))\}
\end{equation}
$Ker(\Psi^{+}_k)$ and $Ker(\Psi^{-}_k)$ are the kernels of the operators $\Psi^{+}_k$ and $\Psi^{-}_k$ ($k$ in $\mathbb{Z}$) as defined in \cite{JPMontillet2013}.  Following Definition $1$, the \emph{uniqueness} of the decomposition  in $\mathbf{s}^{-}(\mathbb{R})$ with the families of differential operators can be stated as:
\vspace{1.0em}
\newline $\bold{Definition}$ $2$: for all $f$ in $\mathbf{s}^{-}(\mathbb{R})$, for all $v\in\mathbb{Z}^+-\{0\}$, for all  $n\in\mathbb{Z}^+$ and $n>1$, the families of operators $(\Psi^{+}_k)_{k \in \mathbb{Z}}$ and $(\Psi^{-}_k)_{k \in \mathbb{Z}}$ ($(\Psi^{+}_k)_{k \in \mathbb{Z}}$ and $(\Psi^{-}_k)_{k \in \mathbb{Z}}$$\subseteq$ $\mathcal{F}(\mathbf{s}^{-}(\mathbb{R}),\mathbf{S}^{-}(\mathbb{R}))$) decompose uniquely $\partial_t^v$ $f^n$ in $\mathbb{R}$, if for any family of operators $(S_k)_{k \in \mathbb{Z}}$ $\subseteq$ $\mathcal{F}(\mathbf{S}^{-}(\mathbb{R}),\mathbf{S}^{-}(\mathbb{R})$) decomposing  $\partial_t^v$$f^n$ in $\mathbb{R}$, there exists a unique couple $(\beta_1,\beta_2)$ in $\mathbb{R}^2$ such as: 
\begin{equation}
S_k(f) = \beta_1 \Psi^{+}_k(f) + \beta_2 \Psi^{-}_k(f), \qquad \forall k\in\mathbb{Z} 
\end{equation}
\vspace{1.0em}
\\ Two important results shown in \cite{JPMontillet2013} are: 
\vspace{1.0em}
\\$\bold{Lemma}$ $0$: for $f$ in $\mathbf{S}^{-}(\mathbb{R})$, the family of DEO ${\Psi}_{k}^{+}$ ($k=\{0,\pm 1,\pm 2,...\}$) decomposes the successive derivatives of the $n$-th power of $f$ for $n\in\mathbb{Z}^+$ and $n>1$. 
\vspace{1.0em}
\\$\bold{Theorem}$ $0$: for $f$ in $\mathbf{s}^{-}(\mathbb{R})$, the families of DEO ${\Psi}_{k}^{+}$ and ${\Psi}_{k}^{-}$ ($k=\{0,\pm 1,\pm 2,...\}$) decompose uniquely the successive derivatives of the $n$-th power of $f$ for $n\in\mathbb{Z}^+$ and $n>1$. 
\vspace{1.0em}
\newline By definition if $f^n$ is analytic, there are ($p$,$q$) ($p >q$) in $\mathbb{R}^2$ such as $f^n$ can be developed in Taylor Series \cite{Kreizig2003}:
\begin{eqnarray}\label{f2eq}
f^n(p) &=& f^n(q) + \sum_{k=1}^\infty \partial_t^k f^n(q) \frac{(p-q)^k}{k!} \nonumber \\
\end{eqnarray}
Let us define for $n$ in $\mathbb{Z}^+-\{0\}$, for $f^n$ in $\mathbf{S}^{-}(\mathbb{R})$ and finite energy, $\mathcal{E}(f^n)$ (in $\mathbf{S}^{-}(\mathbb{R})$) the energy function defined for ($\tau$,$q$) ($q< \tau$) in $\mathbb{R}^2$ such as:
\begin{equation}\label{EnergyfunDefine02}
\mathcal{E}(f^n(\tau)) = \int_q^{\tau} (f^n(t))^2dt < \infty
\end{equation}
\vspace{1.0em}
\newline $\bold{Proposition}$ $1$: If for any $n$ in $\mathbb{Z}^+$, $f^n$ in $\mathbf{S}^{-}(\mathbb{R})$ is analytic and finite energy; for any ($p$,$q$) in $\mathbb{R}^2$ (with $p>q$) and $\mathcal{E}(f^n)$ in $\mathbf{S}^{-}(\mathbb{R})$ is analytic, then $\mathcal{E}(f^n(p))$ is a convergent series.
\begin{proof}
From \eqref{EnergyfunDefine02} and for ($p$,$q$)  in $\mathbb{R}^2$ (with $p>q$), the development in Taylor series of $\mathcal{E}(f^n(p))$ is convergent and can be written as:
\begin{eqnarray}\label{refeq1}
\mathcal{E}(f^n(p)) &=& \mathcal{E}(f^n(q)) + \sum_{k=0}^\infty \partial_t^k (f^n(q))^2 \frac{(p-q)^k}{k!} <\infty \nonumber \\
\end{eqnarray}
Now, let us assume that this series is divergent then from \cite{Kreizig2003},
\begin{eqnarray}
lim_{k\rightarrow \infty} \Big |\frac{\partial_t^k f^{2n}(q)}{\partial_t^{k-1} f^{2n}(q)}\frac{(p-q)}{k+1}\Big| >1
\end{eqnarray}
for $k_1$ and $k_2$ in $\mathbb{Z}^+-\{0\}$ such as $k_1 >>k_2$, 
\begin{eqnarray}
\Big |\frac{\partial_t^{k_1} f^{2n}(q)}{\partial_t^{k_2-1} f^{2n}(q)}\frac{(p-q)^{k_1-k_2}}{(k_1+1)...(k_1-k_2-1)} \Big| &>>&1
\end{eqnarray}
and thus we can conclude clearly  that $\mathcal{E}(f^n(p)) > \infty$. Finally, \eqref{refeq1} is valid if and only if the development in Taylor series of $\mathcal{E}(f^n(p))$ is convergent for $p$ in $\mathbb{R}$.
\end{proof}
%
%
\section{Energy Space}\label{Energyspacesection}
Let us define the open sets $\mathbf{M}^i$ $\subseteq \mathbf{S}^{-}(\mathbb{R})$  for $i$ in $\mathbb{Z}^+$ such as:  
\begin{equation}
\mathbf{M}^i=\{ g \in \mathbf{S}^{-}(\mathbb{R}) | \hspace{0.5em} g = \partial_t^i f^n, f \in \mathbf{s}^{-}(\mathbb{R}), \hspace{0.5em} \hspace{0.5em}  n \in \mathbb{Z}^+-\{0\}\}
\end{equation}
%
Following $\bold{Theorem}$ $0$, $\partial_t^i f^n$ ($i$ in $\mathbb{Z}^+-\{0\}$) can be decomposed uniquely with the family of energy operators $({\Psi}_{k}^{+})_{k\in\mathbb{Z}^+}$ and $({\Psi}_{k}^{-})_{k\in\mathbb{Z}^+}$. It then exists $\alpha_n$ in $\mathbb{R}$ such as:
\begin{equation}
\partial_t^i f^n = \alpha_n ( \partial_t^{i-1} f^{n-2} (\Psi_{1}^{+}(f)+\Psi_{1}^{-}(f)))
\end{equation}
Thus, the above relationship is a linear sum of the energy operators $\Psi_{1}^{+}$ and $\Psi_{1}^{-}$ (and their derivatives)  when $n$ equal $2$. By definition $\Psi_{1}^{-}(f)$ is equal to $0$ for any $f$ in $\mathbf{s}^{-}(\mathbb{R})$. As $\partial_t^i f^n \subseteq \mathcal{F}(\mathbf{s}^{-}(\mathbb{R}),\mathbf{M}^i)$, in the special case $n$ equal $2$, we have $\partial_t^{i-1}\Psi_{1}^{+}(f)$ $\subsetneq \mathcal{F}(\mathbf{s}^{-}(\mathbb{R}),\mathbf{M}^i)$. This is not verified for $n>2$, because the above relationship becomes nonlinear.
\\ One can define an energy space such as:
%
\vspace{1.0em}
\newline $\bold{Definition}$ $3$: The energy space $\mathbf{E}$ is equal to $\mathbf{E}=\bigcup_{i\in\mathbb{Z}^+ \cup \{0\}} \mathbf{M}^i$. 
 \vspace{1.0em}
\\ Note that the definition of $\mathbf{M}^0$ does not involve the family of DEO ${\Psi}_{k}^{+}$ and ${\Psi}_{k}^{-}$ ($k=\{0,\pm 1,\pm 2,...\}$). In this particular case, if $g$ is a general solution of some PDEs, then $f^n$ can be assimilated as some special form of the solution (if it exists). That will be investigated in the last section of this document. 
\\In some specific applications (e.g., solutions of some PDEs), It may be interesting to use the energy function $\mathcal{E}(.)$ applied to $\partial_t^i f^n$ if we want to restrict the definition of $\mathbf{E}$ (not to all $\mathbb{Z}^+$). In this case, $\mathcal{E}(.)$ is said $\emph{associated}$ with the energy space $\mathbf{E}$.
\vspace{1.0em}
\\Furthermore from \cite{Kreizig2003}, the energy of the function $f^n$ is directly connected to the $L1$ norm with the Cauchy-Schwartz inequality (for $(p,q)$ and $q<p$, $\mathbb{R}^2$, $(\int_q^p f^n dt)^2 \leq (p-q)\mathcal{E}(f^n)$).
\vspace{1.0em}
\newline $\bold{Proposition}$ $2$ : for all $f$ in $\mathbf{S}^{-}(\mathbb{R})$, for all $k \in \mathbb{Z}^+$, $\mathcal{E}( \Psi_{k+1}^{+}(f))$ $\leq$ $\mathcal{E}(\partial^k_t \Psi_{k}^{+}(f))$  
\vspace{1.0em}
\begin{proof}
From the properties of the derivative chain rules and with Cauchy-Schwartz inequality \cite{Kreizig2003} one can write :
\begin{eqnarray}
\Psi_{k+1}^{+}(f) &=& \partial_t \Psi_{k}^{+}(f) - \Psi_{k-1}^{+}(\partial_t f) \nonumber \\
\int_{-\infty}^{+\infty} |\Psi_{k+1}^{+}(f)|^2 dt &\leq& \int_{-\infty}^{+\infty} |\partial_t \Psi_{k}^{+}(f)|^2 dt \nonumber \\
\mathcal{E}( \Psi_{k+1}^{+}(f)) &\leq&\mathcal{E}(\partial^k_t \Psi_{k}^{+}(f)) \nonumber \\
\end{eqnarray}
\end{proof}
%
%
\section{Generalized Energy Operators}\label{sectionenergyopdefandgen}
From the introduction of the DEO families in \cite{Maragos1995} and as recalled in Section \ref{Introduction part}, it is possible to generalize the definition of the energy operators $(\Psi_{k}^{+})_{k\in\mathbb{Z}}$ and $(\Psi_{k}^{-})_{k\in\mathbb{Z}}$ based on some operators defined as:
\begin{eqnarray}
{[.,.]}_k^- &=& \partial_t.\partial_t^{k-1}. - . \partial_t^k, \hspace{0.5em} k\in\mathbb{Z} \nonumber \\
{[.,.]}_k^- &=& [.]_k^- \\
\end{eqnarray}
Let us call it the generalized energy operator $[.,.]_k^-$ $\subseteq \mathcal{F}(\mathbf{S}^-(\mathbb{R}),\mathbf{S}^-(\mathbb{R}))$. Note that  in \cite{Boudraa et al.2009}, the authors defined a similar operator using Lie bracket restricted to signal processing applications (e.g., signal and speech AM-FM demodulation). 
Similarly to the DEO families, one can introduce the conjugate $[.,.]_k^+$ $\subseteq \mathcal{F}(\mathbf{S}^-(\mathbb{R}),\mathbf{S}^-(\mathbb{R}))$ defined as:
\begin{eqnarray}\label{bracketl}
{[.,.]}_k^+ &=& \partial_t.\partial_t^{k-1}. + . \partial_t^k., \hspace{0.5em} k\in\mathbb{Z} \nonumber \\
{[.,.]}_k^+ &=& [.]_k^+ \\
\end{eqnarray}
To obtain the families of DEOs $(\Psi_k^+)_{k\in\mathbb{Z}}$ and $(\Psi_k^-)_{k\in\mathbb{Z}}$ defined in \cite{JPMontillet2013}, one can then apply $[.,.]_k^+$ and $[.,.]_k^-$ to $f$ in $\mathbf{S}^{-}(\mathbb{R})$ such as:
\begin{eqnarray}\label{eqpsip2}
{[f,f]}_{k}^+ &=& \partial_t f \partial_t^{k-1} f + f \partial_t^k f \nonumber \\
{[f,f]_k}^+ &=& {\Psi}_k^+(f) \nonumber \\
{[f,f]_k}^- &=& \partial_t f \partial_t^{k-1} f - f \partial_t^k f \nonumber \\
{[f,f]_k}^- &=& {\Psi}_k^-(f) \nonumber \\
\end{eqnarray} 
Furthermore, one can write:
\begin{eqnarray}\label{generalizedEO}
{[[f,f]_k^+,[f,f]_k^+]}_{k}^+ &=& \partial_t {\Psi}_k^+(f) \partial_t^{k-1} {\Psi}_k^+(f) + {\Psi}_k^+(f) \partial_t^k {\Psi}_k^+(f) \nonumber \\
&=&[[f]^1]_{k}^+ \nonumber \\
{[[f,f]_k^-,[f,f]_k^-]}_{k}^- &=& \partial_t {\Psi}_k^-(f) \partial_t^{k-1} {\Psi}_k^-(f) - {\Psi}_k^-(f) \partial_t^k {\Psi}_k^-(f) \nonumber \\
&=&[[f]^1]_{k}^- \nonumber \\
\end{eqnarray} 
Here, we define the notation $[[.]^p]_{k}^+$ with $k$ in $\mathbb{Z}$ and $p$ in $\mathbb{Z}^+$. $k$ is the degree of the derivative similar to the definition of $\Psi_k^+(.)$, and $p$ is the number of "iterations" of the operator $[.,.]_k^+$. Thus following the equation above, $\Psi_k^+(f)$ is equal to $[[f]^0]_{k}^+$, and $\Psi_k^-(f)$ equal to $[[f]^0]_{k}^-$.
%
Now, let us show the proposition:
\vspace{1.0em}
\newline $\bold{Proposition}$ $3$: for all  $p\in\mathbb{Z}^+$, for all $k$ in $\mathbb{Z}$, the generalized energy operators $[[.]^p]_k^+$ and $[[.]^p]_k^-$ are bilinear and follow the derivative chain rule property.
\vspace{1.0em}
\begin{proof}
First, let us recall the properties of a bilinear map according to \cite{Nathan}. Let us define the set $\mathbf{V} \subseteq \mathbf{S}^-(\mathbb{R})$, and the map $B: \mathbf{V} \times \mathbf{V} \rightarrow \mathbf{S}^-(\mathbb{R})$. $B$ is a bilinear map if and only if:
\begin{eqnarray}
B({v}_1 + {v}_2, w) &=& B({v}_1,w) + B({v}_2,w), \hspace{0.5em} \forall {v}_1, {v}_2, w \in \mathbf{V} \nonumber \\
B(v_1, w_1+w_2) &=& B(v_1,w_1) + B(v_1,w_2), \hspace{0.5em} \forall v_1, w_1, w_2 \in \mathbf{V} \nonumber \\
B(cv_1, w) &=& B(v_1, cw) \hspace{0.5em} = \hspace{0.5em} cB(v_1, w), \hspace{0.5em} \forall v_1, w \in \mathbf{V}, \hspace{0.5em} c \in \mathbb{R} \nonumber \\
\end{eqnarray}
\\Previous works such as \cite{Boudraa et al.2009} and \cite{JPMontillet2013}  showed that for any $k$ in $\mathbb{Z}$,  $\Psi_k^+(.)$ and  $\Psi_k^-(.)$  are quadratic forms of a specific bilinear operator. Thus for any $k$ in $\mathbb{Z}$,  $\Psi_k^+(.)$ and  $\Psi_k^-(.)$ are bilinear operators due to the quadratic superposition principle \cite{Boudraa et al.2009}. In addition, with the definition of the family of energy operators $(\Psi_k^+)_{k\in\mathbb{Z}}$ and  $(\Psi_k^-)_{k\in\mathbb{Z}}$, it is straightforward that $([.,.]_k^+)_{k\in\mathbb{Z}}$ and $([.,.]_k^-)_{k\in\mathbb{Z}}$ are families of bilinear operators. In the following, the bilinearity and the derivative chain rule property of the generalized energy operator families are shown by induction on the index $p$ in  $\mathbb{Z}^+$.
\vspace{1.0em}
 \\$A$ - $\bold{Bilinearity}$
\vspace{1.0em}
\begin{itemize}
\item$\bold{Case}$ $p=0$
\end{itemize}
\vspace{1.0em}
 To refer to the previous paragraph, based on $[.,.]_k^+$ ($[.,.]_k^-$) in \eqref{bracketl}, this is the generalization of the quadratic operator $\Psi_k^+(.)$ ($\Psi_k^-(.)$). Therefore, $[[.]^0]_k^+$ and $[[.]^0]_k^-$ are bilinear operators. 
\vspace{1.0em}
\begin{itemize}
\item$\bold{Case}$ $p=1$
\end{itemize}
\vspace{1.0em}
By definition, 
\begin{eqnarray}\label{eqpsip222}
{[ [.]^1 ]}_{k}^+ &=& \big [ [.]_k^0,[.]_k^0\big]_k^+ \nonumber \\
{[ [.]^1 ]}_{k}^+ &=&\big [ \Psi_k^+, \Psi_k^+ \big]_k^+ \nonumber \\
\end{eqnarray} 
With $\Psi_k^+(.)$ ($k$ in $\mathbb{Z}$) and $[.,.]_k^+$ bilinear operators, we can conclude that  $[[.]^1]_{k}^+$ is a bilinear operator as well for any $k$ in $\mathbb{Z}$.
\vspace{1.0em}
\begin{itemize}
\item$\bold{Case}$ $p=h+1$
\end{itemize}
\vspace{1.0em}
Now, we assume that $[[.]^h]_{k}^+$ is a bilinear operator for any $k$ in $\mathbb{Z}$. We can write:
\begin{eqnarray}\label{eqpsip222b}
[[.]^{h+1}]_{k}^+ &=& \big [ [.]_k^{h},[.]_k^{h}\big]_k^+ \nonumber \\
\end{eqnarray} 
As mentioned before, $[[.]^h]_{k}^+$ and $[.,.]_k^+$ are bilinear operators. Thus, we can conclude that  ${[[.]^{h+1}]}_{k}^+$ is also a bilinear operator for any $k$ in $\mathbb{Z}$. By replacing $+$ with $-$ in the previous equations, it shows that the bilinearity of the conjugate operator ${[[.]^p]}_{k}^-$.
\vspace{1.0em}
\\ $B$ - $\bold{Derivative}$ $\bold{Chain}$ $\bold{Rule}$
\vspace{1.0em}
\\ Now, let us show the derivative chain rule property of $[[.]^p]_{k}^-$ and $[[.]^p]_{k}^+$ for any $k$ in $\mathbb{Z}$ and $p$ in $\mathbb{Z}^+$ with induction.
\vspace{1.0em}
\begin{itemize}
\item$\bold{Case}$ $p=0$
\end{itemize}
\vspace{1.0em}
It was shown in \cite{JPMontillet2013} that $[[.]^0]_{k}^+$ and $[[.]^0]_{k}^-$ ($k$ in $\mathbb{Z}$) follow the chain rules derivative rule such as for any $f$ in $\mathbf{S}^{-}(\mathbb{R})$:
\begin{eqnarray}
{[[f]^0]}_{k}^+ &=& [[f]^0]_{k+1}^+ + [\partial_t[f]^0]_{k-1}^+ \nonumber \\
{[[f]^0]}_{k}^- &=& [[f]^0]_{k+1}^- + [\partial_t[f]^0]_{k-1}^- \nonumber \\
\end{eqnarray}
\vspace{1.0em}
\begin{itemize}
\item$\bold{Case}$ $p=1$
\end{itemize}
\vspace{1.0em}
By definition of the generalized energy operator and with \eqref{eqpsip222b}, one can write for any $f$ in $\mathbf{S}^{-}(\mathbb{R})$:
\begin{eqnarray}
 {[[f]^1]}_{k}^+ &=& \big [ \Psi_k^+(f), \Psi_k^+(f) \big ]_k^+ \nonumber \\
  {[[f]^1]}_{k}^+ &=& \partial_t^{k-1} \Psi_k^+(f) \partial_t \Psi_k^+(f) +\Psi_k^+(f) \partial_t^{k} \Psi_k^+(f) \nonumber \\
\partial_t {[[f]^1]}_{k}^+ &=& \partial_t^{k} \Psi_k^+(f) \partial_t \Psi_k^+(f) + \partial_t^{k-1} \Psi_k^+(f) \partial_t^2 \Psi_k^+(f)  +\partial_t \Psi_k^+(f) \partial_t^{k} \Psi_k^+(f) + \Psi_k^+(f) \partial_t^{k+1} \Psi_k^+(f) \nonumber \\
\partial_t {[[f]^1]}_{k}^+ &=& {[[f]^1]}_{k+1}^+ + {[\partial_t[f]^1]}_{k-1}^+ \nonumber \\
\end{eqnarray}
\vspace{1.0em}
\begin{itemize}
\item$\bold{Case}$ $p=h+1$
\end{itemize}
\vspace{1.0em}
Let us assume that the derivative chain rule works for the generalize operator  $[[.]^h]_{k}^+$ ($k$ in $\mathbb{Z}$). Following the previous case, we write:
\begin{eqnarray}
 {[[f]^{h+1}]}_{k}^+ &=& \big [ {[[f]^{h}]}_{k}^+, {[[f]^{h}]}_{k}^+ \big ]_k^+ \nonumber \\
 {[[f]^{h+1}]}_{k}^+ &=& \partial_t^{k-1}  {[[f]^{h}]}_{k}^+ \partial_t  {[[f]^{h}]}_{k}^+ +  {[[f]^{h}]}_{k}^+ \partial_t^{k}  {[[f]^{h}]}_{k}^+ \nonumber \\
 \partial_t  {[[f]^{h+1}]}_{k}^+ &=& \partial_t^{k}  {[[f]^{h}]}_{k}^+ \partial_t  {[[f]^{h}]}_{k}^+ +\partial_t^{k-1}  {[[f]^{h}]}_{k}^+ \partial_t^2  {[[f]^{h}]}_{k}^+ + {[[f]^{h}]}_{k}^+ \partial_t^{k+1}  {[[f]^{h}]}_{k}^+ + \partial_t {[[f]^{h}]}_{k}^+ \partial_t^{k}  {[[f]^{h}]}_{k}^+ \nonumber \\
 \partial_t  {[[f]^{h+1}]}_{k}^+ &=& {[[f]^{h+1}]}_{k+1}^+ + {[\partial_t[f]^{h+1}]}_{k-1}^+ \nonumber \\
\end{eqnarray}
This is the end of the proof by induction of the derivative chain rule for the generalized operator $[[.]^p]_{k}^+$. The same induction can be done for the generalized energy operator $[[.]^p]_{k}^-$ by simply replacing the sign. Note that the derivative chain rule property of the generalized energy operators comes from the general Leibniz derivative rules. This can be compared to similar properties of fractal operators such as mentioned in \cite{Bolognaetal}. 
\end{proof}
Note that for $p$ in $\mathbb{Z}^+$ and $f$ in $\mathbf{S}^{-}(\mathbb{R})$, $[[f]^p]_{k}^+$ in $\mathbf{S}^{-}(\mathbb{R})$, and for $n\in\mathbb{Z}^+$ and $n>1$ $\big ( [[f]^p]_{k}^\pm \big )^n$ in $\mathbf{S}^{-}(\mathbb{R})$. With this property, it is possible to extent $\bold{Lemma}$ $0$ and $\bold{Theorem}$ $0$ (e.g., \cite{JPMontillet2013}) established for the families of DEOs $(\Psi_k^+)_{k\in\mathbb{Z}}$ and $(\Psi_k^-)_{k\in\mathbb{Z}}$ when using the generalized energy operators. 
\vspace{0.5em}
\newline However following the energy space definition in Section \ref{Energyspacesection}, let us introduce, for $p$ in $\mathbb{Z}^+$, the energy space $\mathbf{E}_p = \bigcup_{i\in\mathbb{Z}^+ \cup \{0\}} \mathbf{H}^i$ with $\mathbf{H}^i \subseteq \mathbf{S}^{-}(\mathbb{R})$ such as :
\begin{equation}\label{energyspaceH}
\mathbf{H}^i=\{ g \in \mathbf{S}^{-}(\mathbb{R}), \hspace{0.5em} p\in \mathbb{Z}^+ | \hspace{0.5em} g = \partial_t^i \big ( \big[ [f]^p \big ]_1^+ \big )^n , \big[ [f]^p \big ]_1^+ \in \mathbf{S}^{-}(\mathbb{R}), \hspace{0.5em} n\in\mathbb{Z}^+-\{0\}\}
\end{equation}
One can define for $p$ in $\mathbb{Z}^+$, $\mathbf{S}_p^{-}(\mathbb{R}) \subseteq \mathbf{S}^{-}(\mathbb{R})$, such as:
\begin{equation}
\mathbf{S}_p^{-}(\mathbb{R})=\{ \hspace{0.5em} p\in \mathbb{Z}^+ | \mathbf{E}_p = \bigcup_{i\in\mathbb{Z}^+ \cup \{0\}} \mathbf{H}^i \neq \{0\} \}
\end{equation}
Similar to the comments in Section 3, the energy space $\mathbf{E}_p$ is said  $\emph{associated}$ with $\mathcal{E}( \big[ [.]^p \big ]_1^+)$. Note that Definition $3$ does not involve directly the energy operators. In other words, $\mathbf{E}$ is not equal to $\mathbf{E}_0$.  As an example, the case $p=0$ is defined such as ($i$ in $\mathbb{Z}^+ \cup \{0\}$):
\begin{equation}\label{energyspaceH}
\mathbf{H}^i=\{ g \in \mathbf{S}^{-}(\mathbb{R}), \hspace{0.5em} p\in \mathbb{Z}^+ | \hspace{0.5em} g = \partial_t^i \big ( \Psi_1^+(f) \big )^n , \Psi_1^+(f) \in \mathbf{S}^{-}(\mathbb{R}), \hspace{0.5em} n\in\mathbb{Z}^+-\{0\}\}
\end{equation}
Furthermore, one can define:
\begin{equation}
\mathbf{S}_0^{-}(\mathbb{R})=\{ p=0 | \mathbf{E}_0 = \bigcup_{i\in\mathbb{Z}^+} \mathbf{H}^i \neq \{0\} \}
\end{equation}
Note that for the case $i =0$, it is similar as in the previous case (e.g, definition of $\mathbf{M}^0$) where $\big[ [f]^p \big ]_1^+$ could be considered as a special solution of some PDEs.
\\Here is the Lemma:
\vspace{1.0em}
%
\newline $\bold{Lemma}$ $1$: for $f$ in $\mathbf{S}_{p}^{-}(\mathbb{R})$, $p$ in  $\mathbb{Z}^+$, the families of generalized energy operators $[[.]^{p}]_k^+$ ($k=\{0,\pm 1,\pm 2,...\}$) decompose the successive derivatives of the $n$-th power of $[[f]^{p-1}]_1^+$ for $n\in\mathbb{Z}^+$ and $n>1$.
\begin{proof}
With the convention that ${[[f]^{-1}]}_1^+$ equal $f$, one can see that if $p$ equal $0$, then $\bold{Lemma}$ $1$ is reduced to the case with the families of generalized energy operators ${[[.]^0]}_k^+$ ($k=\{0,\pm 1,\pm 2,...\}$) decomposing the successive derivatives of the $n$-th power of $f$ for $n\in\mathbb{Z}^+$ and $n>1$. This is exactly the statement of $\bold{Lemma}$ $0$. 
Thus,the proof of the $\bold{Lemma}$ $1$ is given by induction on the index $p$ and the the $n$-th power of $[[f]^{p-1}]_1^+$. The induction is used to show the decomposition, and in a separated part on the non-uniqueness. However, this can be long and repetitive compared with the work already published in \cite{JPMontillet2013}. Thus, the $\bold{Lemma}$ $1$ is demonstrated for the case $p=\{0,1,N\}$ ($N$ in $\mathbb{Z}^+$) and $n$ in $\{2,3,L\}$ ($L>1$, $L$ in $\mathbb{Z^+}$).
%
\vspace{1.0em}
\\ $>>$ \emph{A -Decomposition with generalized energy operators} 
\vspace{1.0em}
\vspace{1.0em}
\begin{itemize}
\item$\bold{Case}$ $p=0$
\end{itemize}
\vspace{1.0em}
Following exactly \cite{JPMontillet2013}, the induction on $n$ is separated in two parts: the decomposition with the energy operator families and the non-uniqueness of the decomposition. As all the results are already properly shown in a previous work, we only remind here the main results. Note that $f$ is in $\mathbf{S}_0^{-}(\mathbb{R})$ which according to \cite{JPMontillet2013}, is equal to $\mathbf{S}^{-}(\mathbb{R})$.
\vspace{1.0em}
\\$\bold{Case}$ $n=2$
\vspace{1.0em}
\newline We showed that when $n=2$, one can decompose $\partial_t^s f^2$ ($s$ in $\mathbb{Z}^+$, $s>0$) with the energy operators $[[.]^{0}]_k^+$ ($k=\{0,\pm 1,\pm 2,...\}$) following the Equation (15) in \cite{JPMontillet2013}:
\begin{eqnarray}\label{coefficientsAPgeneral}
\partial_t^s f^2 &=& a_s^+(f) \nonumber \\
\partial_t^s f^2 &=& \sum_{k=0}^{s-1} \big(_{k}^{s-1} \big) {[\partial_t^{s-k-1}[f]^{0}]}_{2(k+1)-s}^{+}, \qquad \forall s\in \mathbb{Z}^{+}-\{0\} 
\end{eqnarray}
\vspace{1.0em}
\\$\bold{Case}$ $n=3$
\vspace{1.0em}
\\ From the Equation ($20$) in \cite{JPMontillet2013} and \eqref{coefficientsAPgeneral}, one can write:
\begin{eqnarray}\label{Newrec01}
\partial_t f^3 = f \frac{3}{2} a_1^+(f) \nonumber \\
\partial_t f^3 = f \frac{3}{2}[[f]^{0}]_1^+  \nonumber \\
\partial_t f^3 = f  A_1^{+}(f) \nonumber \\
& & \nonumber \\
A_s^+(f) = \frac{3}{2} \partial_t^{s-1} a_1^+(f), \qquad \forall s\in \mathbb{Z}^{+}-\{0\} \nonumber \\
A_s^+(f) = \frac{3}{2} \partial_t^{s-1} [[f]^{0}]_1^+, \qquad \forall s\in \mathbb{Z}^{+}-\{0\} \nonumber \\
& & \nonumber \\
\partial_t^{s} f^3 = \sum_{k=0}^{s-1} \big(_{k}^{s-1} \big) A_{k+1}^+(f) \partial_t^{s-1-k} f, \qquad \forall s\in \mathbb{Z}^+-\{0\} 
\end{eqnarray}
\vspace{1.0em}
\\$\bold{Case}$ $n=L$, $L>1$
\vspace{1.0em}
\\ With the notation of the generalized energy operators, as shown in Equation ($28$) in \cite{JPMontillet2013}:
\begin{eqnarray}\label{NONOVO}
 \partial_t f^L =  p f^{L-1} \partial_t f  \nonumber \\
 \partial_t f^L =  \frac{L}{2} f^{L-2} [[f]^{0}]_1^+  \nonumber \\
\partial_t f^L = \frac{L}{L-1} B_1^+(f)  f^{L-2}\nonumber \\
& & \nonumber \\
B_s^+(f) = \frac{L-1}{2}\partial_t^{s-1} [[f]^{0}]_1^+, \qquad \forall s\in \mathbb{Z}^{+}-\{0\} \nonumber \\
& & \nonumber \\
%
\partial_t^2 f^L = \frac{L}{L-1} B_1^+(f) \partial_t f^{L-2}+ \frac{L}{L-1} B_2^+(f)  f^{L-2} \nonumber \\
\partial_t^{s+1} f^L =\sum_{k=0}^s \big(_{k}^s \big) \frac{L}{L-1}B_{k+1}^+(f) \partial_t^{s-k} f^{L-2}, \qquad \forall s \in \mathbb{Z}^+
\end{eqnarray}
That is all the main results shown in the case of $\bold{Lemma}$ $0$. This then ends the case $p=0$. 
\vspace{1.0em}
\begin{itemize}
\item$\bold{Case}$ $p=1$
\end{itemize}
\vspace{1.0em}
In this case, the $\bold{Lemma}$ $1$ should be stated as:  for $f$ in $\mathbf{S}_1^{-}(\mathbb{R})$,  the families of generalized energy operators $[[.]^{1}]_k^+$ ($k=\{0,\pm 1,\pm 2,...\}$) decompose the successive derivatives of the $n$-th power of $[[f]^{0}]_1^+$ for $n\in\mathbb{Z}^+$ and $n>1$.
\newline In other words, one can substitute $f^n$ with $\big ( [[f]^{0}]_1^+ \big )^n$, and the family of generalized operator $\big ( [[.]^{0}]_k^+ \big )_{k\in\mathbb{Z}^+}$ with  $\big ( [[.]^{1}]_k^+ \big )_{k\in\mathbb{Z}^+}$. It is then possible to write according to the previous development:
\vspace{1.0em}
\\$\bold{Case}$ $n=2$
\vspace{1.0em}
\newline Following the same development as in the previous case, 
\begin{eqnarray}\label{coefficientsAPgeneralvv}
\partial_t^s \big ( [[f]^{0}]_1^+ \big )^2 &=& a_s^+(f) \nonumber \\
\partial_t^s \big ( [[f]^{0}]_1^+ \big )^2 &=& \sum_{k=0}^{s-1} \big(_{k}^{s-1} \big) {[\partial_t^{s-k-1}[f]^{1}]}_{2(k+1)-s}^{+}, \qquad \forall s\in \mathbb{Z}^{+}-\{0\} 
\end{eqnarray}
\vspace{1.0em}
\\$\bold{Case}$ $n=3$
\vspace{1.0em}
\\ From the  \eqref{Newrec01}, one can write:
\begin{eqnarray}
\partial_t \big ( [[f]^{0}]_1^+ \big )^3 = f \frac{3}{2}[[f]^{1}]_1^+  \nonumber \\
\partial_t \big ( [[f]^{0}]_1^+ \big )^3 = f \frac{3}{2} a_1^+(f) \nonumber \\
\partial_t \big ( [[f]^{0}]_1^+ \big )^3 = f  A_1^{+}(f) \nonumber \\
& & \nonumber \\
A_s^+(f) = \frac{3}{2} \partial_t^{s-1} a_1^+(f), \qquad \forall s\in \mathbb{Z}^{+}-\{0\} \nonumber \\
A_s^+(f) = \frac{3}{2} \partial_t^{s-1} [[f]^{1}]_1^+, \qquad \forall s\in \mathbb{Z}^{+}-\{0\} \nonumber \\
& & \nonumber \\
\partial_t^{s} \big ( [[f]^{0}]_1^+ \big )^ 3 = \sum_{k=0}^{s-1} \big(_{k}^{s-1} \big) A_{k+1}^+(f) \partial_t^{s-1-k} f, \qquad \forall s\in \mathbb{Z}^+-\{0\} 
\end{eqnarray}
\vspace{1.0em}
\\$\bold{Case}$ $n=L$, $L>1$
\vspace{1.0em}
\\ With the notation of the generalized energy operators, it was shown in \eqref{NONOVO}:
\begin{eqnarray}\label{NONOVOvv}
 \partial_t \big ( [[f]^{0}]_1^+ \big )^L =  L \big ( [[f]^{0}]_1^+ \big )^{L-1} \partial_t \big ( [[f]^{0}]_1^+ \big )  \nonumber \\
 \partial_t \big ( [[f]^{0}]_1^+ \big )^L =  \frac{L}{2} \big ( [[f]^{0}]_1^+ \big )^{L-2} [[f]^{1}]_1^+  \nonumber \\
\partial_t \big ( [[f]^{0}]_1^+ \big )^L = \frac{L}{L-1} B_1^+(f)  \big ( [[f]^{0}]_1^+ \big )^{L-2}\nonumber \\
& & \nonumber \\
B_s^+(f) = \frac{L-1}{2}\partial_t^{s-1} [[f]^{1}]_1^+, \qquad \forall s\in \mathbb{Z}^{+}-\{0\} \nonumber \\
& & \nonumber \\
%
\partial_t^2 \big ( [[f]^{0}]_1^+ \big )^L = \frac{L}{L-1} B_1^+(f) \partial_t \big ( [[f]^{0}]_1^+ \big )^{L-2}+ \frac{L}{L-1} B_2^+(f)  \big ( [[f]^{0}]_1^+ \big )^{L-2} \nonumber \\
\partial_t^{s+1} \big ( [[f]^{0}]_1^+ \big )^L =\sum_{k=0}^s \big(_{k}^s \big) \frac{L}{L-1}B_{k+1}^+(f) \partial_t^{s-k} \big ( [[f]^{0}]_1^+ \big )^{L-2}, \qquad \forall s \in \mathbb{Z}^+
\end{eqnarray}
\vspace{1.0em}
\begin{itemize}
\item$\bold{Case}$ $p=N$
\end{itemize}
\vspace{1.0em}
In this case, the $\bold{Lemma}$ $1$ states that:  for $f$ in $\mathbf{S}_N^{-}(\mathbb{R})$,  the families of generalized energy operators $[[.]^{N}]_k^+$ ($k=\{0,\pm 1,\pm 2,...\}$) decompose the successive derivatives of the $n$-th power of $[[f]^{N-1}]_1^+$ for $n\in\mathbb{Z}^+$ and $n>1$.
\\ As we see in the statement of the $Lemma$ $1$, one has to assume that $\mathbf{S}_N^{-}(\mathbb{R})$ is not reduced to $\{0\}$. Thus, following the previous development:
\vspace{1.0em}
\\$\bold{Case}$ $n=2$
\vspace{1.0em}
\newline Following the same development as in the previous case, 
\begin{eqnarray}\label{coefficientsAPgeneralvv2}
\partial_t^s \big ( [[f]^{N-1}]_1^+ \big )^2 &=& a_s^+(f) \nonumber \\
\partial_t^s \big ( [[f]^{N-1}]_1^+ \big )^2 &=& \sum_{k=0}^{s-1} \big(_{k}^{s-1} \big) {[\partial_t^{s-k-1}[f]^{N}]}_{2(k+1)-s}^{+}, \qquad \forall s\in \mathbb{Z}^{+}-\{0\} \nonumber \\
\end{eqnarray}
\vspace{1.0em}
\\$\bold{Case}$ $n=3$
\vspace{1.0em}
\\ From \eqref{Newrec01}, one can write:
\begin{eqnarray}
\partial_t \big ( [[f]^{N-1}]_1^+ \big )^3 = f \frac{3}{2}[[f]^{N}]_1^+  \nonumber \\
\partial_t \big ( [[f]^{N-1}]_1^+ \big )^3 = f \frac{3}{2} a_1^+(f) \nonumber \\
\partial_t \big ( [[f]^{N-1}]_1^+ \big )^3 = f  A_1^{+}(f) \nonumber \\
& & \nonumber \\
A_s^+(f) = \frac{3}{2} \partial_t^{s-1} a_1^+(f), \qquad \forall s\in \mathbb{Z}^{+}-\{0\} \nonumber \\
A_s^+(f) = \frac{3}{2} \partial_t^{s-1} [[f]^{N}]_1^+, \qquad \forall s\in \mathbb{Z}^{+}-\{0\} \nonumber \\
& & \nonumber \\
\partial_t^{s} \big ( [[f]^{N-1}]_1^+ \big )^ 3 = \sum_{k=0}^{s-1} \big(_{k}^{s-1} \big) A_{k+1}^+(f) \partial_t^{s-1-k} f, \qquad \forall s\in \mathbb{Z}^+-\{0\} 
\end{eqnarray}
\vspace{1.0em}
\\$\bold{Case}$ $n=L$, $L>1$
\vspace{1.0em}
With the notation of the generalized energy operators, it was shown in \eqref{NONOVO}:
\begin{eqnarray}\label{NONOVOvv2}
 \partial_t \big ( [[f]^{N-1}]_1^+ \big )^L =  L \big ( [[f]^{N-1}]_1^+ \big )^{L-1} \partial_t \big ( [[f]^{N-1}]_1^+ \big )  \nonumber \\
 \partial_t \big ( [[f]^{N-1}]_1^+ \big )^L =  \frac{L}{2} \big ( [[f]^{N-1}]_1^+ \big )^{L-2} [[f]^{N}]_1^+  \nonumber \\
\partial_t \big ( [[f]^{N-1}]_1^+ \big )^L = \frac{L}{L-1} B_1^+(f)  \big ( [[f]^{N-1}]_1^+ \big )^{L-2}\nonumber \\
& & \nonumber \\
B_s^+(f) = \frac{L-1}{2}\partial_t^{s-1} [[f]^{N}]_1^+, \qquad \forall s\in \mathbb{Z}^{+}-\{0\} \nonumber \\
& & \nonumber \\
%
\partial_t^2 \big ( [[f]^{N-1}]_1^+ \big )^L = \frac{L}{L-1} B_1^+(f) \partial_t \big ( [[f]^{N-1}]_1^+ \big )^{L-2}+ \frac{L}{L-1} B_2^+(f)  \big ( [[f]^{N-1}]_1^+ \big )^{L-2} \nonumber \\
\partial_t^{s+1} \big ( [[f]^{N-1}]_1^+ \big )^L =\sum_{k=0}^s \big(_{k}^s \big) \frac{L}{L-1}B_{k+1}^+(f) \partial_t^{s-k} \big ( [[f]^{N-1}]_1^+ \big )^{L-2}, \qquad \forall s \in \mathbb{Z}^+
\end{eqnarray}
This ends the first part on the decomposition of functions with generalized energy operators.
\vspace{1.0em}
\\ $>>$ \emph{B - Non-uniqueness of the decomposition with generalized energy operators} 
\vspace{1.0em}
%
\begin{itemize}
\item$\bold{Case}$ $p=0$
\end{itemize}
\vspace{1.0em}
In \cite{JPMontillet2013} in the proof of the $Lemma$ $0$, the non-uniqueness of the decomposition of the successive derivatives $\partial_t^i f^n$ ($n\in \mathbb{Z}^+$, $n>1$, $i\in\mathbb{Z}^+$)  was shown with a simple counter example for $f$ in  $\mathbf{S}_0^-(\mathbb{R})$:
\begin{eqnarray}\label{definitionOfeta}
 \eta_k(f) &=& 3(\partial_t f \partial_t^{k-1} f)-f\partial_t^k f , \hspace{0.5em} \forall k \in \mathbb{Z} \nonumber \\
\end{eqnarray}
Note that the derivative chain rule property is applied to this operator. One can verify:
\begin{eqnarray}\label{Coeff2nnn}
\partial_t f^2 &=& 2 f \partial_t f \nonumber \\
\partial_t f^2 &=& \eta_1(f) \nonumber \\ 
\eta_1(f)&=&{[[f]^0]}_{1}^{+} \nonumber \\
& & \nonumber \\
\partial_t^2 f^2 &=& 2 (\partial_t f)^2 + 2 f \partial_t^2 f \nonumber \\
\partial_t^2 f^2 &=&  \partial_t {\eta}_{1}^{+}(f) \nonumber \\
\partial_t {[[f]^0]}_{1}^{+}  &=& \partial_t {\eta}_{1}^{+}(f) \nonumber \\
\partial_t^2 f^2 &=& {\eta}_{2}^{+}(f) + {\eta}_{0}^{+}(\partial_t f) \nonumber \\
\end{eqnarray}
\vspace{1.0em}
\begin{itemize}
\item$\bold{Case}$ $p=1$
\end{itemize}
\vspace{1.0em}
In the same way, one can also define the generalized energy operator $[[.]^{1}]_{\eta k}^+$ ($k=\{0,\pm 1,\pm 2,...\}$) decomposing the successive derivatives of the $n$-th power of $[[f]^{0}]_1^+$  for $f$ in  $\mathbf{S}_1^-(\mathbb{R})$, for $n\in\mathbb{Z}^+$ and $n>1$.
\begin{eqnarray}\label{Coeff2nnnetak}
[[f]^{1}]_{\eta k}^+ &=& 3(\partial_t [[f]^{0}]_1^+ \partial_t^{k-1} [[f]^{0}]_1^+)- [[f]^{0}]_1^+\partial_t^k [[f]^{0}]_1^+ , \qquad \forall k \in \mathbb{Z} \nonumber \\
& & \nonumber \\
\partial_t \big ( [[f]^{0}]_1^+ \big )^2 &=& 2 [[f]^{0}]_1^+ \partial_t f \nonumber \\
\partial_t \big ( [[f]^{0}]_1^+ \big )^2 &=& [[f]^{1}]_{\eta 1}^+ \nonumber \\ 
{[[f]^{1}]}_{\eta 1}^+ &=&[[f]^{1}]_{1}^+ \nonumber \\
& & \nonumber \\
\partial_t^2 \big ( [[f]^{0}]_1^+ \big )^2 &=& 2 (\partial_t [[f]^{0}]_1^+ )^2 + 2 [[f]^{0}]_1^+ \partial_t^2 [[f]^{0}]_1^+ \nonumber \\
\partial_t^2 \big ( [[f]^{0}]_1^+ \big )^2 &=&  \partial_t [[f]^{1}]_{\eta 1}^+ \nonumber \\
\partial_t [[f]^{1}]_{1}^+ &=& \partial_t [[f]^{1}]_{\eta 1}^+ \nonumber \\
\partial_t^2 \big ( [[f]^{0}]_1^+ \big )^2 &=& [[f]^{1}]_{\eta 2}^+ + [\partial_t[f]^{1}]_{\eta 0}^+ \nonumber \\
\end{eqnarray}
\vspace{1.0em}
\begin{itemize}
\item$\bold{Case}$ $p=N$
\end{itemize}
\vspace{1.0em}
One can generalize the case $p=\{0,1\}$ with the generalized energy operator  for $f$ in  $\mathbf{S}_N^-(\mathbb{R})$: 
\begin{equation}\label{Coeff2nnnetakb2}
{[[f]^{N}]}_{\eta k}^+ = 3(\partial_t [[f]^{N-1}]_1^+ \partial_t^{k-1} [[f]^{N-1}]_1^+)- [[f]^{N-1}]_1^+\partial_t^k [[f]^{N-1}]_1^+ , \qquad \forall k \in \mathbb{Z} \nonumber \\
\end{equation}
Following the same development,
\begin{eqnarray}\label{Coeff2nnnetakN}
\partial_t \big ( [[f]^{N-1}]_1^+ \big )^2 &=& 2 [[f]^{N-1}]_1^+ \partial_t f \nonumber \\
\partial_t \big ( [[f]^{N-1}]_1^+ \big )^2 &=& [[f]^{N}]_{\eta 1}^+ \nonumber \\ 
{[[f]^{N}]}_{\eta 1}^+ &=&[[f]^{N}]_{1}^+ \nonumber \\
& & \nonumber \\
\partial_t^2 \big ( [[f]^{N-1}]_1^+ \big )^2 &=& 2 (\partial_t [[f]^{N-1}]_1^+ )^2 + 2 [[f]^{N-1}]_1^+ \partial_t^2 [[f]^{N-1}]_1^+ \nonumber \\
\partial_t^2 \big ( [[f]^{N-1}]_1^+ \big )^2 &=&  \partial_t [[f]^{N}]_{\eta 1}^+ \nonumber \\
\partial_t [[f]^{N}]_{1}^+ &=& \partial_t [[f]^{N}]_{\eta 1}^+ \nonumber \\
\partial_t^2 \big ( [[f]^{N-1}]_1^+ \big )^2 &=& [[f]^{N}]_{\eta 2}^+ + [\partial_t[f]^{N}]_{\eta 0}^+ \nonumber \\
\end{eqnarray}
\end{proof}
Now, one can define the subset $\mathbf{s}_p^{-}(\mathbb{R})$ defined as:
\begin{equation}
\mathbf{s}_p^{-}(\mathbb{R}) = \{ f \in  \mathbf{S}_p^{-}(\mathbb{R}), \hspace{0.5em} p \in \mathbb{Z}^+ | f \notin  (\cup_{k \in \mathbb{Z}} Ker([[f]^{p}]_k^+)\cup(\cup_{k \in \mathbb{Z}-\{1\}} Ker([[f]^{p}]_k^-))\}
\end{equation}
The subset $\mathbf{s}_p^{-}(\mathbb{R})$ is also defined such as $\mathbf{E}_p\neq \{0\}$. Thus, one can see that $\mathbf{s}_p^{-}(\mathbb{R})\subseteq\mathbf{S}_p^{-}(\mathbb{R})$.
\vspace{1em}\newline $\bold{Theorem}$ $1$: for $f$ in $\mathbf{s}_p^{-}(\mathbb{R})$, for $p$ in $\mathbb{Z}^+$, the families of generalized operators $[[.]^p]_{k}^+$ and $[[.]^p]_{k}^-$ ($k=\{0,\pm 1,\pm 2,...\}$) decompose uniquely the successive derivatives of the $n$-th power of $[[f]^{p-1}]_{1}^+$ for $n\in\mathbb{Z}^+$ and $n>1$. 
\begin{proof}
The proof is an induction on the index $p$ and the $n$-th power of $[[f]^{p-1}]_{1}^+$. It is separated in three parts: the decomposition with the generalized energy operators, the existence and the uniqueness of the decomposition. However, that follows exactly the work of \cite{JPMontillet2013}. Similarly to the proof of $\bold{Lemma}$ $1$, one should notice that for the case $p=0$ the $\bold{Theorem}$ $1$ is exactly the statement of $\bold{Theorem}$ $0$ with $\mathbf{s}_0^{-}(\mathbb{R})$ equal to $\mathbf{s}^{-}(\mathbb{R})$ defined in the first section. To keep the demonstration short, the induction is done for $n$ in $\{2,L\}$ and $p$ in $\{0,N\}$. 
\vspace{1.0em}
\\ $>>$ \emph{A -Decomposition with generalized energy operators} 
\vspace{1.0em}
\vspace{1.0em}
\begin{itemize}
\item$\bold{Case}$ $p=0$
\end{itemize}
\vspace{1.0em}
Recall the definition of $[[.]^0]_{1}^-$ from \eqref{generalizedEO} and the proof of $\bold{Theorem}$ $0$ (see \cite{JPMontillet2013}), one can write for $f$ in $\mathbf{s}_0^-(\mathbb{R})$:
\vspace{1.0em}
\newline $\bold{Case}$ $n=2$:
\vspace{1.0em}
\begin{eqnarray}\label{equationPsikm01}
\partial_t f^2 &=& f \partial_t f + f\partial_t f + f\partial_t f - f\partial_t f \nonumber \\
\partial_t f^2 &=&[[f]^0]_{1}^+ + [[f]^0]_{1}^-\nonumber \\
& & \nonumber \\
\partial_t^2 f^2 &=& 2 (\partial_t f)^2 + 2 f \partial_t^2 f \nonumber \\
\partial_t^2 f^2 &=& \partial_t ([[f]^0]_{1}^+ +[[f]^0]_{1}^-) \nonumber \\
\partial_t^2 f^2 &=& [[f]^0]_{2}^++ [\partial_t [f]^0]_{0}^++ \nonumber \\
                 &&  [[f]^0]_{2}^-+ [\partial_t [f]^0]_{0}^-\nonumber \\
\end{eqnarray}
There is a symmetry with the proof of the previous lemma. From \eqref{equationPsikm01}, one can define $a_s^-(f)$ in the same way that $a_s^+(f)$ was defined in \eqref{coefficientsAPgeneral} as:
\begin{equation}\label{Aocorf1N3}
\partial_t^s f^2 = \sum_{k=0}^{s-1} \big(_{k}^{s-1} \big) {[\partial_t^{s-k-1}{[f]}^{0}]}_{2(k+1)-s}^{+}, \qquad \forall s\in \mathbb{Z}^{+}-\{0\}
\end{equation}
%
%
%
With the property of the derivative chain rule, it is easy to calculate the first terms of the DEO family $a_s^-(f)$ such as :
\begin{eqnarray}\label{Coeff2nmoins}
%
a_1^-(f)=[[f]^0]_{1}^- &=& 0 \nonumber \\
& & \nonumber \\
a_2^-(f)=\partial_t[[f]^0]_{1}^- &=&[[f]^0]_{2}^- + [\partial_t [f]^0]_{0}^- \nonumber \\
{[[f]^0]}_{2}^- &=& -{[\partial_t {[f]}^0]}_{0}^- \nonumber \\
& & \nonumber \\
\end{eqnarray}
The family of DEO $[[f]^0]_{0}^-$ ($k\in\mathbb{Z}$) has  the same derivative properties as $[[f]^0]_{0}^+$. A similar equation can then be established for $a_s^-(f)$ following the development written in \eqref{coefficientsAPgeneral} as:
\begin{eqnarray}\label{coefficientsAPgeneralbis}
  a_s^-(f)&=& \sum_{k=0}^{s-1} \big(_{k}^{s-1} \big) {[\partial_t^{s-k-1}[f]^{0}]}_{2(k+1)-s}^{-}, \qquad \forall s\in \mathbb{Z}^{+}-\{0\}  \nonumber \\
  a_s^-(f) &=& 0 \nonumber \\
\end{eqnarray}
This formula has just been checked for $s=\{1,2,3,4\}$ with \eqref{Coeff2nmoins}.
The generalization of the formula for $s=m$ is very similar to the one described in \eqref{coefficientsAPgeneralbis} literally by changing $+$ and $-$ in the definition of the energy operator. It follows that the decomposition of the successive derivatives of $f^2$ is generalized for any $n$ in $\mathbb{Z}^+-\{0\}$ as:
\begin{eqnarray}\label{Coeff2nmoinspp}
\partial_t^m f^2 &=& \partial_t^{m-1} ([[f]^0]_{0}^+ +[[f]^0]_{0}^-) \nonumber \\
                 &=& \sum_{k=0}^{m-1} \big(_{k}^{m-1} \big) {[\partial_t^{m-k-1}[f]^{0}]}_{2(k+1)-m}^{+} + \nonumber \\
                 & & \sum_{k=0}^{m-1} \big(_{k}^{m-1} \big) {[\partial_t^{m-k-1}[f]^{0}]}_{2(k+1)-m}^{-} \nonumber \\
\end{eqnarray}
\vspace{1.0em}
\\$\bold{Case}$ $n=L$, $L>1$
\vspace{1.0em}
\\ Following the same step as in the proof of $\bold{Lemma}$ $1$, let us define $B_s^+(f)$ and $B_s^-(f)$ ( $s$ in $\mathbb{Z}^{+}-\{0\} $) with the assumption that they decompose the successive derivatives of ${f}^{p-1}$ as:
\begin{eqnarray}\label{psi+0008}
\partial_t f^L &=&  p f^{L-1} \partial_t f  \nonumber \\
\partial_t f^L &=&  \frac{L}{2} f^{L-2} [[f]^{0}]_1^+  \nonumber \\
%
%
& & \nonumber \\
B_s^+(f) &=& \frac{L-1}{2}\partial_t^{s-1} [[f]^{0}]_1^+, \qquad \forall s\in \mathbb{Z}^{+}-\{0\} \nonumber \\
B_s^-(f) &=& \frac{L-1}{2}\partial_t^{s-1} [[f]^{0}]_1^-, \qquad \forall s\in \mathbb{Z}^{+}-\{0\} \nonumber \\
%
\partial_t f^L &=& \frac{L}{L-1} \big (B_1^+(f)+B_1^-(f)\big )  f^{L-2}\nonumber \\
& & \nonumber \\
%
%
\partial_t^2 f^L &=& \frac{L}{L-1} \big (B_1^+(f)+B_1^-(f)\big ) \partial_t f^{L-2}+ \frac{L}{L-1} \big (B_2^+(f)+B_2^-(f)\big )  f^{L-2} \nonumber \\
%
%
\end{eqnarray}
There is again a symmetry with the proof of $\bold{Lemma}$ $1$.  One can define the $s+1$-th derivative of $\big ( [[f]^{N-1}]_1^+ \big )^L$ using $B_{k+1}^-(f)$ and $B_{k+1}^+(f)$: 
\begin{equation}\label{dtppfwithApsd2}
\partial_t^{s+1} f^L =\sum_{k=0}^s \big(_{k}^s \big) \frac{L}{L-1} \big (B_{k+1}^+(f)+ B_{k+1}^-(f) \big ) \partial_t^{s-k} f^{L-2}, \qquad \forall s \in \mathbb{Z}^+
\end{equation}
This equation has just been checked for $s=\{0,1\}$. As the induction proof follows exactly the proof of $\bold{Lemma}$ $1$ as in \eqref{NONOVO} by only adding $B_{k+1}^-(f)$ which has the same properties as $B_{k+1}^+(f)$. It allows then to assume the generalization to the case $s+2$.
\newline Thus,  $(B_{k}^+)_{k\in\mathbb{Z}}$ and  $(B_{k}^-)_{k\in\mathbb{Z}}$ decompose the $s$-th derivative of $f^L$. From their definition, one can conclude that $\big ( [[f]^{0}]_k^+ \big )_{k\in\mathbb{Z}}$ and $\big ( [[f]^{0}]_k^- \big )_{k\in\mathbb{Z}}$ decompose $\partial_t^s f^L$.
\vspace{1.0em}
\begin{itemize}
\item$\bold{Case}$ $p=N$
\end{itemize}
\vspace{1.0em}
%
This case follows the proof of the $\bold{Lemma}$ $1$ and in particular \eqref{coefficientsAPgeneralvv2} to \eqref{NONOVOvv2} and \eqref{Coeff2nmoins}. One can then write for $f$ in $\mathbf{s}_N^-(\mathbb{R})$:
\vspace{1.0em}
\newline $\bold{Case}$ $n=2$:
\vspace{1.0em}
\begin{eqnarray}\label{coefficientsAPgeneralvv3}
\partial_t^s \big ( [[f]^{N-1}]_1^+ \big )^2 &=& a_s^+(f) + a_s^-(f)\nonumber \\
\partial_t^s \big ( [[f]^{N-1}]_1^+ \big )^2 &=& \sum_{k=0}^{s-1} \big(_{k}^{s-1} \big) \big ( {[\partial_t^{s-k-1}[f]^{N}]}_{2(k+1)-s}^{+} +{[\partial_t^{s-k-1}[f]^{N}]}_{2(k+1)-s}^{-}\big ), \qquad \forall s\in \mathbb{Z}^{+}-\{0\} \nonumber \\
\end{eqnarray}
\vspace{1.0em}
\vspace{1.0em}
\\$\bold{Case}$ $n=L$, $L>1$
\vspace{1.0em}
With the notation of the generalized energy operators, it was shown in \eqref{NONOVOvv}:
\begin{eqnarray}\label{NONOVOvv3}
 \partial_t \big ( [[f]^{N-1}]_1^+ \big )^L &=&  L \big ( [[f]^{N-1}]_1^+ \big )^{L-1} \partial_t \big ( [[f]^{N-1}]_1^+ \big )  \nonumber \\
 \partial_t \big ( [[f]^{N-1}]_1^+ \big )^L &=&  \frac{L}{2} \big ( [[f]^{N-1}]_1^+ \big )^{L-2} \big ([[f]^{N}]_1^+  +[[f]^{N}]_1^- \big ) \nonumber \\
\partial_t \big ( [[f]^{N-1}]_1^+ \big )^L &=& \frac{L}{L-1} (B_1^+(f)+ B_1^-(f)) \big ( [[f]^{N-1}]_1^+ \big )^{L-2}\nonumber \\
& & \nonumber \\
B_s^+(f) &=& \frac{L-1}{2}\partial_t^{s-1} [[f]^{N}]_1^+, \qquad \forall s\in \mathbb{Z}^{+}-\{0\} \nonumber \\
B_s^-(f) &=& \frac{L-1}{2}\partial_t^{s-1} [[f]^{N}]_1^-, \qquad \forall s\in \mathbb{Z}^{+}-\{0\} \nonumber \\
& & \nonumber \\
\partial_t^2 \big ( [[f]^{N-1}]_1^+ \big )^L &=& \frac{L}{L-1} (B_1^+(f)+ B_1^-(f))  \partial_t \big ( [[f]^{N-1}]_1^+ \big )^{L-2} \nonumber \\
& &+ \frac{L}{L-1} (B_2^+(f)+ B_2^-(f))  \big ( [[f]^{N-1}]_1^+ \big )^{L-2} \nonumber \\
\partial_t^{s+1} \big ( [[f]^{N-1}]_1^+ \big )^L &=&\sum_{k=0}^s \big(_{k}^s \big) \frac{L}{L-1}(B_{k+1}^+(f)+ B_{k+1}^-(f))  \partial_t^{s-k} \big ( [[f]^{N-1}]_1^+ \big )^{L-2}, \qquad \forall s \in \mathbb{Z}^+\nonumber \\
\end{eqnarray}
\vspace{1.0em}
\\ $>>$ \emph{B -Existence of the decomposition with generalized energy operators} 
\vspace{1.0em}
\\ In the proof of $\bold{Lemma}$ $1$, it was shown that the non-uniqueness of the decomposition using a counter-example. Here, the proof re-investigate these examples. Then, it is generalized for $p >0$ via induction on $p$.  
\vspace{1.0em}
\begin{itemize}
\item$\bold{Case}$ $p=0$
\end{itemize}
\vspace{1.0em}
With $f$ in $\mathbf{s}_0^-(\mathbb{R})$,
\\$\bold{Case}$ $n=2$: 
\vspace{1.0em}
\\It was shown that the family of operators $(\eta_k)_{k\in\mathbb{Z}}$ (proof of the Lemma, \eqref{definitionOfeta}), decomposes $\partial^s_t f^{2}$ ($s\in\mathbb{Z}^+-\{0\}$).  One can rewrite it as a sum of the DEO family $({\Psi}_{k}^{-})_{k\in\mathbb{Z}}$ and $({\Psi}_{k}^{+})_{k\in\mathbb{Z}}$ as:
\begin{eqnarray}
\eta_k(f) &=& \Psi_k^+(f) + 2 \Psi_k^-(f), \qquad k \in \mathbb{Z} \nonumber \\
\eta_k(f) &=&[[f]^{0}]_k^+ +2 [[f]^{0}]_k^- , \qquad k \in \mathbb{Z} \nonumber \\
\end {eqnarray}
\vspace{1.0em}
\\$\bold{Case}$ $n=L$ with $L>1$: 
%
%
\\Previously, \eqref{psi+0008} defined the decomposition of $\partial^{s}_t f^{L}$ ($s\in\mathbb{Z}^+-\{0\}$) with the generalized energy operator $[[.]^{0}]_k^+$. With the definition of $[[.]^{0}]_1^-$, one can define the operator $\theta_k^+$ and $\theta_k^-$ as:
\begin{eqnarray}\label{psi+000100}
\partial_t f^L &=& L  f^{L-1} \partial_t f\nonumber \\
\partial_t f^L &=&  \frac{L}{2} f^{L-2} [[f]^{0}]_1^+  \nonumber \\
& & \nonumber \\
\theta_k^+(f) &=&  \frac{L-1}{2}  [[f]^{0}]_k^+ \nonumber \\
\theta_k^-(f) &=&  \frac{L-1}{2}  [[f]^{0}]_k^- \nonumber \\
& & \nonumber \\
\partial_t f^L &=&  \frac{L}{2} f^{L-2} \big ( [[f]^{0}]_1^+ +[[f]^{0}]_1^- \big ) \nonumber \\
&=& \frac{L}{L-1} f^{L-2}(\theta_1^+(f) + \theta_1^-(f))
\end{eqnarray}
Following the development in \eqref{psi+0008} and \eqref{NONOVOvv3}, one can see that $B_s^+(f) = \partial_t^{s-1} \theta_1^+(f)$ and $B_s^-(f) = \partial_t^{s-1} \theta_1^-(f)$ ($s\in\mathbb{Z}^+-\{0\}$). Note that ${\theta}_{k}^{-}$ and ${\theta}_{k}^{+}$ are bilinear operators and follow the derivative chain rule property by definition.
\newline Using  \eqref{NONOVOvv3}, one can easily show that ${\theta}_{k}^{-}$ and ${\theta}_{k}^{+}$ decomposes $\partial^{s}_t f^{L}$ ($s\in\mathbb{Z}^+-\{0\}$). It is then possible to conclude the existence of the decomposition of any operator by using $([[.]^{0}]_k^-)_{k\in\mathbb{Z}}$ and $([[.]^{0}]_k^+)_{k\in\mathbb{Z}}$.
\vspace{1.0em}
\begin{itemize}
\item$\bold{Case}$ $p=N$
\end{itemize}
\vspace{1.0em}
With $f$ in $\mathbf{s}_N^-(\mathbb{R})$,
\\$\bold{Case}$ $n=2$: 
\vspace{1.0em}
\\It was shown that the family of operators $({[[.]^{N}]}_{\eta k}^+)_{k\in\mathbb{Z}}$ (proof of $\bold{Lemma}$ $1$, \eqref{Coeff2nnnetakb2}), decomposes $\partial^s_t ([[f]^{N-1}]_1^+)^{2}$ ($s\in\mathbb{Z}^+-\{0\}$).  One can rewrite it as a sum of the DEO family $([[.]^{N}]^-_{k})_{k\in\mathbb{Z}}$ and $([[.]^{N}]^+_{k})_{k\in\mathbb{Z}}$ as:
\begin{eqnarray}
{[[f]^{N}]}_{\eta k}^+ &=& 3(\partial_t [[f]^{N-1}]_1^+ \partial_t^{k-1} [[f]^{N-1}]_1^+)- [[f]^{N-1}]_1^+\partial_t^k [[f]^{N-1}]_1^+ \nonumber \\
{[[f]^{N}]}_{\eta k}^+ &=& [[f]^{N}]_1^+ + 2 [[f]^{N}]_1^- , \qquad k \in \mathbb{Z} \nonumber \\
\end {eqnarray}
\vspace{1.0em}
\\$\bold{Case}$ $n=L$ with $L>1$: 
%
%
\\Previously, \eqref{psi+0008} defined the decomposition of $\partial^{s}_t \big( [[.]^{N-1}]_1^+\big )^{L}$ ($s\in\mathbb{Z}^+-\{0\}$) with the generalized energy operator $[[.]^{N}]_k^+$. With the definition of $[[.]^{N}]_1^-$, one can define the operator $[[.]^{N}]_{\theta k}^+$ and $[[.]^{N}]_{\theta k}^-$ as:
\begin{eqnarray}\label{psi+000100bc}
\partial_t \big( [[f]^{N-1}]_1^+\big )^{L} &=& L  \big( [[f]^{N-1}]_1^+\big )^{L-1} \partial_t \big( [[f]^{N-1}]_1^+\big )\nonumber \\
\partial_t \big( [[f]^{N-1}]_1^+\big )^{L} &=&  \frac{L}{2} \big( [[f]^{N-1}]_1^+\big )^{L-2}  [[f]^{N}]_1^+  \nonumber \\
& & \nonumber \\
{[[f]^{N}]}_{\theta k}^+ &=&  \frac{L-1}{2}  {[[f]^{N}]}_k^+ \nonumber \\
{[[f]^{N}]}_{\theta k}^- &=&  \frac{L-1}{2}  {[[f]^{N}]}_k^- \nonumber \\
& & \nonumber \\
\partial_t \big( [[f]^{N-1}]_1^+\big )^{L}  &=& \frac{L}{L-1} f^{L-2} f^{L-2} \big ( [[f]^{N}]_{\theta 1}^+ +[[f]^{N}]_{\theta 1}^- \big ) \nonumber \\
\end{eqnarray}
Following the development in \eqref{psi+0008} and \eqref{NONOVOvv3}, one can see that $B_s^+(f) = \partial_t^{s-1} [[f]^{N}]_{\theta 1}^+$ and $B_s^-(f) = \partial_t^{s-1} [[f]^{N}]_{\theta 1}^-$ ($s\in\mathbb{Z}^+-\{0\}$). Note that ${[[f]^{N}]}_{\theta k}^-$ and ${[[f]^{N}]}_{\theta k}^+$ are bilinear operators and follow the derivative chain rule property by definition.
\newline With \eqref{psi+000100bc}, one can easily show that ${[[.]^{N}]}_{\theta k}^-$ and ${[[.]^{N}]}_{\theta k}^+$ decomposes  $\partial^{s}_t \big( [[.]^{N-1}]_1^+\big )^{L}$ ($s\in\mathbb{Z}^+-\{0\}$). It is then possible to conclude the existence of the decomposition of any operator by using $([[.]^{N}]_k^-)_{k\in\mathbb{Z}}$ and $([[.]^{N}]_k^+)_{k\in\mathbb{Z}}$.
\vspace{1.0em}
\\ $>>$ \emph{C - About the Uniqueness of the decomposition with generalized energy operators} 
\vspace{1.0em}
\\ Following the previous sections, the proof by induction on the index $p$ in $\mathbb{Z}^+$ shows the uniqueness of the decomposition of any family of operators decomposing $\partial_t^{s+1} \big ([[f]^{p-1}]_k^+\big)^n$ ($f$ in $\mathbf{s}_p^-(\mathbb{R})$, $s$ in $\mathbb{Z}^+$, $n$ in $\mathbb{Z}^+$ and $n > 1$) with the families of generalized operators $[[.]^p]_{k}^+$ and $[[.]^p]_{k}^-$ ($k=\{0,\pm 1,\pm 2,...\}$) . In other words for example for $p$ is equal to $0$, one wants to show that if a family of operators $({S}_{k})_{k\in\mathbb{Z}}$ ( $S_k \subsetneq \mathcal{F}( \mathbf{s}_0^{-}(\mathbb{R}), \mathbf{S}^{-}(\mathbb{R}))$) decomposes $\partial_t^{s+1} f^n$ ($s$ in $\mathbb{Z}^+$, $n$ in $\mathbb{Z}^+$ and $n > 1$), ${S}_{k}$ ($k$ in $\mathbb{Z}$) can be written with an unique sum of $[[.]^0]_{k}^+$ and $[[.]^0]_{k}^-$ ($k=\{0,\pm 1,\pm 2,...\}$). Thus, the induction is on the index $p$ and the $k$-th order of the generalized energy operators. Note that for a matter of clarity, $p$ is restricted to the case $\{0,N\}$.
\vspace{1.0em}
\begin{itemize}
\item $\bold{Case}$ $p=0$
\end{itemize}
\vspace{1.0em}
This is the case already shown in the proof of $\bold{Theorem}$ $0$. The same logic of the proof is applied for the case $p>0$.
\vspace{1.0em}
\\$\bold{Case}$ $k=2$: 
For $f$ in $\mathbf{s}_0^{-}(\mathbb{R})$ (or  $\mathbf{s}^{-}(\mathbb{R})$) and $n$ in $\mathbb{Z}^+$ and $n > 1$, one can assume that $(\alpha_1,\alpha_2,\beta_1,\beta_2)$ exist in $\mathbb{R}^4$ such as:
%
%
%
%
\begin{eqnarray}\label{eqrefpartial0102}
\partial_t^s f^n &=& \partial_t^{s-1} S_1(f)  \nonumber \\
\partial_t^s f^n &=& \alpha_1 \partial_t^{s-1} [[f]^0]_{1}^+ +\alpha_2 \partial_t^{s-1} [[f]^0]_{1}^- \nonumber \\
\partial_t^s f^n &=& \beta_1 \partial_t^{s-1} [[f]^0]_{1}^+ +\beta_2 \partial_t^{s-1} [[f]^0]_{1}^- \nonumber \\
\end{eqnarray} 
%
%
\\ As  with the operator family $({S}_{k})_{k\in\mathbb{Z}}$ follows the derivative chain rule property:
\begin{eqnarray}\label{equationalphabeta1}
\partial_t S_1(f) &=& S_2(f) +S_0(\partial_t f)\nonumber \\
\partial_t S_1(f) &=&\alpha_1 \partial_t [[f]^0]_{1}^+ +\alpha_2 \partial_t [[f]^0]_{1}^- \nonumber \\
\partial_t S_1(f) &=&\alpha_1 ([[f]^0]_{2}^+ +[\partial_t[f]^0]_{0}^+) +\alpha_2 ([[f]^0]_{2}^- +[\partial_t[f]^0]_{0}^-) \nonumber \\
\end{eqnarray} 
And then,
\begin{eqnarray}\label{equationalphabeta1b}
S_2(f) &=& \alpha_1 [[f]^0]_{2}^+ +\alpha_2 [[f]^0]_{2}^- \nonumber \\
S_2(f) &=& \beta_1 [[f]^0]_{2}^+ +\beta_2 [[f]^0]_{2}^- \nonumber \\
(\alpha_1-\beta_1) [[f]^0]_{2}^+ +(\alpha_2-\beta_2) [[f]^0]_{2}^- &=& 0
\end{eqnarray} 
As $f$ in $\mathbf{s}_0^{-}(\mathbb{R})$,  the images of $[[.]^0]_{2}^+$  and $[[.]^0]_{2}^-$ ($Im([[.]^0]_{2}^+)$ and $Im([[.]^0]_{2}^-)$) are not reduced to $\{0\}$. It follows that $\alpha_1=\beta_1$ and $\alpha_2=\beta_2$. Note that it is not possible to do this simple check for $k=1$ as $Im([[.]^0]_{1}^-)=\{0\}$ by definition of the family  $\big ( [[.]^p]_{k}^- \big )_{k \in \mathbb{Z}}$ .
\vspace{1.0em}
\\$\bold{Case}$ $k=L$: 
For $f$ in $\mathbf{s}_0^{-}(\mathbb{R})$, let us assume the uniqueness of the decomposition for $k=L-1$ (with $k\neq1$). For $k=L$, following \eqref{equationalphabeta1}:
\begin{eqnarray}\label{equationalphabeta1kn}
\partial_t S_{L-1}(f) &=& S_L(f) +S_{L-2}(\partial_t f)\nonumber \\
\partial_t S_{L-1}(f) &=&\alpha_1 \partial_t \ [[f]^0]_{L-1}^+ +\alpha_2 \partial_t  [[f]^0]_{L-1}^+ \nonumber \\
\partial_t S_{L-1}(f) &=&\alpha_1 ( [[f]^0]_{L}^+ + [\partial_t[f]^0]_{L-2}^+)+\alpha_2 ( [[f]^0]_{L}^- + [\partial_t[f]^0]_{L-2}^-) \nonumber \\
\end{eqnarray} 
And then,
\begin{eqnarray}\label{equationalphabeta1knb}
S_L(f) &=& \alpha_1 [[f]^0]_{L}^+ +\alpha_2 [[f]^0]_{L}^- \nonumber \\
S_L(f) &=& \beta_1 [[f]^0]_{L}^+ +\beta_2 [[f]^0]_{L}^- \nonumber \\
(\alpha_1-\beta_1) [[f]^0]_{L}^++(\alpha_2-\beta_2) [[f]^0]_{L}^- &=& 0
\end{eqnarray} 
By definition for $L\neq1$, $Im([[f]^0]_{L}^+)$ and $Im( [[f]^0]_{L}^-)$ are not reduced to $\{0\}$, and it follows that $\alpha_1=\beta_1$ and $\alpha_2=\beta_2$. 
\vspace{1.0em}
\\$\bold{Special}$ $\bold{Case}$ $k=1$: 
To complete the proof with the assumption that  $\alpha_1=\beta_1$ and $\alpha_2=\beta_2$ for $k \in\mathbb{Z}$ and $k \neq 1$ , the special case $k=1$ can be solved as:
\begin{eqnarray}\label{equationalphabeta1c}
\partial_t(\alpha_1 [[f]^0]_{1}^+ ) &=&\alpha_1([[f]^0]_{2}^+ + [\partial_t[f]^0]_{0}^+) \nonumber \\
&=&\beta_1([[f]^0]_{2}^+ + [\partial_t[f]^0]_{0}^+) \nonumber \\
&=&\partial_t(\beta_1[[f]^0]_{1}^+)\nonumber \\
& & \nonumber \\
\partial_t(\alpha_2 [[f]^0]_{1}^-) &=&\alpha_2([[f]^0]_{2}^- + [\partial_t[f]^0]_{0}^-)\nonumber \\
&=&\beta_2([[f]^0]_{2}^- + [\partial_t[f]^0]_{0}^-)\nonumber \\
&=&\partial_t(\beta_2 [[f]^0]_{1}^-)\nonumber \\
\nonumber\\
\end{eqnarray} 
To conclude in equation ($53$) in \cite{JPMontillet2013}, it is shown that $\alpha_1 = \beta_1 =1$.
\vspace{1.0em}
\begin{itemize}
\item$\bold{Case}$ $p=N$
\end{itemize}
\vspace{1.0em}
In this case, $f$ is in $\mathbf{s}_N^{-}(\mathbb{R})$. Following the previous development, one can assume that there is a family of energy operators $(V_k)_{k\in\mathbb{Z}}$ and $V_k \subsetneq \mathcal{F}( \mathbf{s}_N^{-}(\mathbb{R}), \mathbf{S}^{-}(\mathbb{R}))$ which decomposes $( [[f]^{N-1}]_1^+\big)^n$. 
\vspace{1.0em}
\\$\bold{Case}$ $k=2$: 
For  $n$ in $\mathbb{Z}^+$ and $n > 1$, one can assume that $(\alpha_1,\alpha_2,\beta_1,\beta_2)$ exist in $\mathbb{R}^4$ such as:
%
%
%
%
\begin{eqnarray}\label{eqrefpartial0102bv}
\partial_t^s \big ( [[f]^{N-1}]_1^+\big)^n &=& \partial_t^{s-1} V_1(f)  \nonumber \\
\partial_t^s f^n &=& \alpha_1 \partial_t^{s-1} [[f]^{N}]_1^+ +\alpha_2 \partial_t^{s-1} [[f]^{N}]_1^- \nonumber \\
\partial_t^s f^n &=& \beta_1 \partial_t^{s-1} [[f]^{N}]_1^+ +\beta_2 \partial_t^{s-1} [[f]^{N}]_1^- \nonumber \\
\end{eqnarray} 
%
%
\\ In addition, the family $({V}_{k})_{k\in\mathbb{Z}}$ follows the derivative chain rule property:
\begin{eqnarray}\label{equationalphabeta1bv}
\partial_t V_1(f) &=& V_2(f) +V_0(\partial_t f)\nonumber \\
\partial_t V_1(f) &=&\alpha_1 \partial_t [[f]^{N}]_1^+ +\alpha_2 \partial_t [[f]^{N}]_1^- \nonumber \\
\partial_t V_1(f) &=&\alpha_1 ([[f]^{N}]_2^+ +[\partial_t[f]^{N}]_0^+ ) +\alpha_2 ([[f]^{N}]_2^- +[\partial_t[f]^{N}]_0^- ) \nonumber \\
\end{eqnarray} 
And then,
\begin{eqnarray}\label{equationalphabeta1bbv}
V_2(f) &=& \alpha_1 [[f]^N]_{2}^+ +\alpha_2 [[f]^N]_{2}^- \nonumber \\
V_2(f) &=& \beta_1 [[f]^N]_{2}^+ +\beta_2 [[f]^N]_{2}^- \nonumber \\
(\alpha_1-\beta_1) [[f]^N]_{2}^+ +(\alpha_2-\beta_2) [[f]^N]_{2}^- &=& 0
\end{eqnarray} 
As $f$ in $\mathbf{s}_N^{-}(\mathbb{R})$,  the images of $[[.]^N]_{2}^+$  and $[[.]^N]_{2}^-$ are not reduced to $\{0\}$. It follows that $\alpha_1=\beta_1$ and $\alpha_2=\beta_2$. Note that it is not possible to do this simple check for $k=1$ as $Im([[.]^N]_{1}^-)=\{0\}$ by definition .
\vspace{1.0em}
\\$\bold{Case}$ $k=L$: 
For $f$ in $\mathbf{s}_N^{-}(\mathbb{R})$, let us assume the uniqueness of the decomposition for $k=L-1$ (with $k\neq1$). For $k=L$, following \eqref{equationalphabeta1}:
\begin{eqnarray}\label{equationalphabeta1knbv}
\partial_t V_{L-1}(f) &=& V_L(f) +V_{L-2}(\partial_t f)\nonumber \\
\partial_t V_{L-1}(f) &=&\alpha_1 \partial_t \ [[f]^N]_{L-1}^+ +\alpha_2 \partial_t  [[f]^N]_{L-1}^+ \nonumber \\
\partial_t V_{L-1}(f) &=&\alpha_1 ( [[f]^N]_{L}^+ + [\partial_t[f]^N]_{L-2}^+)+\alpha_2 ( [[f]^N]_{L}^- + [\partial_t[f]^N]_{L-2}^-) \nonumber \\
\end{eqnarray} 
And then,
\begin{eqnarray}\label{equationalphabeta1knbbv}
V_L(f) &=& \alpha_1 [[f]^N]_{L}^+ +\alpha_2 [[f]^N]_{L}^- \nonumber \\
V_L(f) &=& \beta_1 [[f]^N]_{L}^+ +\beta_2 [[f]^N]_{L}^- \nonumber \\
(\alpha_1-\beta_1) [[f]^N]_{L}^++(\alpha_2-\beta_2) [[f]^N]_{L}^- &=& 0
\end{eqnarray} 
By definition for $L\neq1$, $Im([[f]^N]_{L}^+)$ and $Im( [[f]^N]_{L}^-)$ are not reduced to $\{0\}$, and it follows that $\alpha_1=\beta_1$ and $\alpha_2=\beta_2$. 
\vspace{1.0em}
\\$\bold{Special}$ $\bold{Case}$ $k=1$: 
To complete the proof with the assumption that  $\alpha_1=\beta_1$ and $\alpha_2=\beta_2$ for $k \in\mathbb{Z}$ and $k \neq 1$ , the special case $k=1$ can be solved as:
\begin{eqnarray}\label{equationalphabeta1c}
\partial_t(\alpha_1 [[f]^N]_{1}^+ ) &=&\alpha_1([[f]^N]_{2}^+ + [\partial_t[f]^N]_{0}^+) \nonumber \\
&=&\beta_1([[f]^N]_{2}^+ + [\partial_t[f]^N]_{0}^+) \nonumber \\
&=&\partial_t(\beta_1[[f]^N]_{1}^+)\nonumber \\
& & \nonumber \\
\partial_t(\alpha_2 [[f]^N]_{1}^-) &=&\alpha_2([[f]^N]_{2}^- + [\partial_t[f]^N]_{0}^-)\nonumber \\
&=&\beta_2([[f]^N]_{2}^- + [\partial_t[f]^N]_{0}^-)\nonumber \\
&=&\partial_t(\beta_2 [[f]^N]_{1}^-)\nonumber \\
\nonumber\\
\end{eqnarray} 
This concludes the proof of $\bold{Theorem}$ $1$.
\end{proof}
%
%
$\bold{Discussion}$ $n<-1$: In this case, one can define:
%
\begin{equation}\label{setbbb}
\forall f \in \mathbf{S}_p^{-}(\mathbb{R}), \hspace{0.5em} \forall t \in \mathbb{R}, \hspace{0.5em} p \in \mathbb{Z}^+, \hspace{0.5em} ([[f(t)]^p]_1^+\big )^n\neq 0, \hspace{0.5em} \forall n \in \mathbb{Z}^+, n>1, \frac{1}{\big ([[f(t)]^p]_1^+\big )^n}
\end{equation}
This set of functions can also be described as: $f$ in $\mathbf{S}_p^{-}(\mathbb{R})$ and $f$ not in  $Ker\big ([[f(t)]^p]_1^+\big )$ for $p$ in $\mathbb{Z}^+$. Note that one could also chose to have $f$ in $\mathbf{s}_p^{-}(\mathbb{R})$. However, this is more restrictive than the set defined in \eqref{setbbb}. Using an intermediary function, $h$ such as $h = \frac{1}{[[f(t)]^p]_1^+}$, the problem of decomposing $\partial_t^s \big ([[f(t)]^p]_1^+\big )^{-n}$ ($s$ in $\mathbb{Z}^+-\{0\}$) is equivalent to resolving $\partial_t^s h^{n}$, which has been demonstrated in the $\bold{Lemma}$ $1$ and $\bold{Theorem}$ $1$.
\vspace{1.0em}
\newline $\bold{Discussion}$ $n=1$ or $n=-1$: As already underlined in \cite{JPMontillet2013}, one can use a general formula for $f$ in the set defined in \eqref{setbbb}:
%
%
\begin{eqnarray}\label{discussion2a}
\partial_t^s \big ([[f(t)]^p]_1^+\big ) &=& \partial_t^s \bigg (\frac{\big ([[f(t)]^p]_1^+\big )^3}{\big ([[f(t)]^p]_1^+\big )^2}\bigg) \nonumber \\
s=1, \qquad \partial_t \big ([[f(t)]^p]_1^+\big ) &=& \big ([[f(t)]^p]_1^+\big )^{-2} \partial_t \big ([[f(t)]^p]_1^+\big )^3 + \big ([[f(t)]^p]_1^+\big )^3 \partial_t \big ([[f(t)]^p]_1^+\big )^{-2} \nonumber \\
s=2, \qquad \partial_t^2 \big ([[f(t)]^p]_1^+\big ) &=& 2 \partial_t \big ([[f(t)]^p]_1^+\big )^{-2} \partial_t \big ([[f(t)]^p]_1^+\big )^3 \nonumber \\
&& + \big ([[f(t)]^p]_1^+\big )^3 \partial_t^2 \big ([[f(t)]^p]_1^+\big )^{-2} + \big ([[f(t)]^p]_1^+\big )^{-2} \partial_t^2 \big ([[f(t)]^p]_1^+\big )^3 \nonumber \\
\end{eqnarray} 
The example for $s =\{1,2\}$ in \eqref{discussion2a} shows that $\partial_t^s \big ([[f(t)]^p]_1^+\big )$  can be decomposed into a product of successive derivatives of $\big ([[f(t)]^p]_1^+\big )^{3}$ and $ \big ([[f(t)]^p]_1^+\big )^{-2}$. Those derivatives can be decomposed into a sum of generalized energy operators based on the $\bold{Lemma}$ $1$ and $\bold{Theorem}$ $1$ plus the previous discussion (for the case $n<-1$ ).
\\ Now for the case $n=-1$, it is easy to see that :
\begin{eqnarray}\label{discussion2amoins1}
\partial_t^s \big ([[f(t)]^p]_1^+\big )^{-1} &=& \partial_t^s \bigg (\frac{\big ([[f(t)]^p]_1^+\big )^2}{\big ([[f(t)]^p]_1^+\big )^3}\bigg) \nonumber \\
\end{eqnarray}
With the discussion for the case $n=1$, we can conclude that $\partial_t^s \big ([[f(t)]^p]_1^+\big )^{-1}$  can be decomposed into a product of successive derivatives of $\big ([[f(t)]^p]_1^+\big )^{2}$ and $ \big ([[f(t)]^p]_1^+\big )^{-3}$.
\section{Solutions of linear PDEs using the energy operators }\label{section4}
In this section and the remainder of this work, the finite energy functions of one variable described in Section \ref{preliminariesSection} (e.g., \eqref{SRRRRR}), are now functions of two variables referring to  the $1$ space dimension ($x$) and time ($t$). Thus, one has to add in the notation of the operators the symbol $t$ or $x$ to indicate which variable the derivatives refer to. For example,  the operators $\Psi_k^{-,x}(.)$ and $[[.]^p]_k^{+,x}$ ($k$ in $\mathbb{Z}$, $p$ in $\mathbb{Z}^+$) refer to their derivatives in space, whereas $\Psi_k^{-,t}(.)$ and $[[.]^p]_k^{+,t}$ ($k$ in $\mathbb{Z}$, $p$ in $\mathbb{Z}^+$) refer to their derivatives in time. This notation agrees with the work in \cite{JPMontillet2010}.
One can then define the Schwartz space $\mathbf{S}^{-}(\mathbb{R}^2)$ for function of two variables such as:
\begin{equation}\label{multipleSR}
\begin{split}
\mathbf{S}^{-}(\mathbb{R}^2) =\{f \in \mathbf{C}^{\infty}(\mathbb{R}),  \hspace{0.5em}  \forall (x_0, t_0) \in \mathbb{R}^+| \hspace{0.5em}  {sup}_{t<0} |t^k||\partial_t^j f(x_0, t)|<\infty, \\
& \hspace{-24em} and \hspace{0.5em} {sup}_{x<0} |x^k||\partial_x^j f(x, t_0)|<\infty ,\hspace{0.5em} \forall k \in \mathbb{Z}^+, \hspace{0.5em} \forall j \in \mathbb{Z}^+ \}
\end{split}
\end{equation}
Following this definition, the extension of the subspace $\mathbf{s}^{-}(\mathbb{R}^2) \subseteq \mathbf{S}^{-}(\mathbb{R}^2)$ is:
\begin{eqnarray}
\mathbf{s}^{-}(\mathbb{R}^2) &=& \{f\in \mathbf{S}^{-}(\mathbb{R}^2)|  \hspace{0.2em}\forall\hspace{0.2em} k \in \mathbb{Z}, \hspace{0.2em}  \Psi^{+,t}_k(f) \neq \{0\} \hspace{0.2em}  \nonumber \\
                   & & \hspace{0.2em} \forall\hspace{0.2em} k \in \mathbb{Z}-\{1\}, \hspace{0.2em}\Psi^{-,t}_k(f) \neq \{0\} \} \bigcup \nonumber  \\
                   & & \{f\in \mathbf{S}^{-}(\mathbb{R}^2)|  \hspace{0.2em}\forall\hspace{0.2em} k \in \mathbb{Z}, \hspace{0.2em}  \Psi^{+,x}_k(f) \neq \{0\} \hspace{0.2em}  \nonumber \\
                   & & \hspace{0.2em} \forall\hspace{0.2em} k \in \mathbb{Z}-\{1\}, \hspace{0.2em}\Psi^{-,x}_k(f) \neq \{0\} \}  \nonumber
\end{eqnarray}
We can redefine Definition $0$, Definition $1$ and the statements of Lemma $0$, Lemma $1$, Theorem $0$ and Theorem $1$ with the function of two variables using the above definitions. We will not state formally all the work previously done for the case of  the functions of two variables in  $\mathbf{S}^{-}(\mathbb{R}^2)$. It is only a matter of replacing the variables from time to space. Note that in some case, if the notations $(\Psi^{+}_k)_{k \in \mathbb{Z}}$ and $(\Psi^{-}_k)_{k \in \mathbb{Z}}$ are used, it means that we are dealing with the families of energy operators in time and in space without making any difference.
\\ With Theorem $1$ in the case of a function of two variables,  it is possible to write the derivatives in space and time with the families of energy operators $(\Psi^{+}_k)_{k \in \mathbb{Z}}$ and $(\Psi^{-}_k)_{k \in \mathbb{Z}}$. One can state that it exists $({\alpha_n}_1, {\alpha_n}_2)$ in $\mathbb{R}^2$ such as for $ f \in \mathbf{s}^{-}(\mathbb{R}^2), \hspace{0.5em} n \in \mathbb{Z}^+ -\{0\}$:
\begin{equation}\label{Problemsetup2}
 \left\{
\begin{array}{rl}
\partial_t^i f^n = & {\alpha_n}_1 ( \partial_t^{i-1} f^{n-2} (\Psi_{1}^{t,+}(f)+\Psi_{1}^{t,-}(f))), \\
\partial_x^i f^n = & {\alpha_n}_2 ( \partial_t^{i-1} f^{n-2} (\Psi_{1}^{x,+}(f)+\Psi_{1}^{x,-}(f))) \\
\end{array} \right.
\end{equation}
%
%
%
In addition, we did not use the generalized energy operator notation (see Section \ref{sectionenergyopdefandgen} ) for the case of the functions of two variables in order to keep it readable. Note that n the same way $\mathbf{S}^{-}(\mathbb{R}^2) $ in \eqref{multipleSR}  is defined for two variables, it can be extended to functions of multiple variables with the linearity of the definition. 
 \\ The remainder of this section shows the definition of the sets of solutions using energy operators and generalized energy operators for linear PDEs.
\subsection{Application to linear PDEs with the families of DEOs}\label{sectionAppPDEs}
Let us consider the linear partial differential equation :
\begin{equation}\label{ProblemsetupB}
 \left\{
\begin{array}{rl}
 a_1 \partial_x^{\beta}g(x,t) + a_2 \partial_t^{\beta}g(x,t) &=h(x,t) ,\\
\beta \in \mathbb{Z}^+-\{0\},& \hspace{0.5em} h \in C^{\infty}(\mathbb{R}^2), \hspace{0.5em} (a_1, a_2) \in \mathbb{R}^2, \\ 
 x\in \mathbb{R}, \hspace{0.5em}  x_0 \in \mathbb{R}, &  t\in \mathbb{R}^+,   \hspace{0.5em} t_0 \in \mathbb{R}^+\\
\end{array} \right.
\end{equation}
$g$ is the general solution of \eqref{ProblemsetupB}. Note that $t$ is now in $\mathbb{R}^+$ as it is more intuitive to define solutions in this interval.
Let us define the open subset $\mathbf{X}^i$ and $\mathbf{Y}^i$  $ \subsetneq  \mathbf{S}^{-}(\mathbb{R}^2)$ with $i$ in $\mathbb{Z}^+$:
\begin{equation}
\mathbf{X}^i =\{  g \in \mathbf{S}^{-}(\mathbb{R}^2), \forall t \in \mathbb{R}^+, \hspace{0.5em} x_0 \in \mathbb{R} | g(x_0,t) = \partial_t^i u(x_0,t)^n,  \hspace{0.5em} u  \in \mathbf{S}^{-}(\mathbb{R}^2),  \hspace{0.5em}n \in \mathbb{Z}^+ -\{0\} \}
\end{equation} 
and,
\begin{equation}
\mathbf{Y}^i =\{  g \in \mathbf{S}^{-}(\mathbb{R}^2), \forall x \in \mathbb{R}, \hspace{0.5em} t_0 \in \mathbb{R}^+ | g(x,t_0) = \partial_x^i u(x_0,t)^n,  \hspace{0.5em} u  \in \mathbf{S}^{-}(\mathbb{R}^2), \hspace{0.5em} n \in \mathbb{Z}^+ -\{0\} \}
\end{equation} 
Note that $\mathbf{M}^0 \subseteq \mathbf{X}^0$ following the discussion after the Definition $3$ (p. 5). The set of all solutions is then defined as $\mathcal{S}_1(\mathbb{R}^2 )= \big ( \bigcup_{i \in \mathbb{Z}^+}  \mathbf{X}^i \big ) \bigcup \big ( \bigcup_{i \in \mathbb{Z}^+}  \mathbf{Y}^i \big )$. With this definition, one can state that $\mathcal{S}_1(\mathbb{R}^2 ) \subsetneq  \mathbf{S}^{-}(\mathbb{R}^2)$. Following the definition of the energy space $\mathbf{E}$ in Definition $3$, $\mathcal{S}_1(\mathbb{R}^2 )$ can also be called an energy space for functions of two variables.
\newline \indent It is important to underline that solutions in $\mathbf{S}^{-}(\mathbb{R}^2)$ of \eqref{ProblemsetupB} are finite energy functions such as the ones decaying for large values of $x$. This is a very limiting condition as we cannot include solutions such as planar waves. 
%
%
\newline \indent Looking at the solutions of \eqref{ProblemsetupB} in   $\mathcal{S}_1(\mathbb{R}^2 )$ can get another meaning when writing the Taylor series expansion of the solution $u$ (in $\mathbf{S}^{-}(\mathbb{R}^2)$)  in time and in space such as for $i$ in $\mathbb{Z}^+$,  $x_0$ in $\mathbb{R}$, and $t_0$ in $\mathbb{R}^+$ :
\begin{equation*}\label{ProblemsetupU}
 \left\{
\begin{array}{rl}
 u^n(x,t_0) = & u^n(x_0,t_0) + \sum_{i=1}^\infty \partial_x^i u^n(x_0,t_0) \frac{(x-x_0)^i}{i!},\\
 u^n(x_0,t) = & u^n(x_0,t_0) + \sum_{i=1}^\infty \partial_t^i u^n(x_0,t_0) \frac{(t-t_0)^i}{i!} \\
\end{array} \right.
\end{equation*}
By definition, one can write $ \partial_t^i u^n(x,t_0) \in \mathbf{X}^i$ and $ \partial_x^i u^n(x_0,t) \in \mathbf{Y}^i$.
\newline \indent In addition if $u$ in  $\mathbf{s}^{-}(\mathbb{R}^2)$,  the derivatives of $u^n$  can be written with the family of energy operators $(\Psi_k^{+,t})_{k\in\mathbb{Z}}$, $(\Psi_k^{-,t})_{k\in\mathbb{Z}}$, $(\Psi_k^{+,x})_{k\in\mathbb{Z}}$ and $(\Psi_k^{-,x})_{k\in\mathbb{Z}}$ following $\bold{Theorem}$ $0$ and the development of the proof of $\bold{Theorem}$ $1$ (see equations \eqref{Coeff2nmoinspp} and \eqref{dtppfwithApsd2}). With the development in the previous section, one can write ($x_0$ in $\mathbb{R}$, and $t_0$ in $\mathbb{R}^+$):
\begin{eqnarray}\label{dtppfwithApsdg2}
\partial_t^{\beta} u^n(x_0,t) &=& \sum_{k=0}^{\beta-1} \big(_{k}^{\beta-1} \big) \frac{n}{2}(\partial_t^k\Psi_1^{+,t}(u)(x_0,t) +\partial_t^k\Psi_1^{-,t}(u))(x_0,t) \nonumber \\
& & \partial_t^{\beta-1-k} u^{n-2}(x_0,t), \hspace{0.5em} \forall \beta \in \mathbb{Z}^+ -\{0\}, \hspace{0.5em} n>1 \\
\partial_x^{\beta} u^n(x,t_0) &=& \sum_{k=0}^{\beta-1} \big(_{k}^{\beta-1} \big) \frac{n}{2}(\partial_x^k\Psi_1^{+,x}(u)(x,t_0) +\partial_x^k\Psi_1^{-,x}(u))(x,t_0) \nonumber \\
& & \partial_x^{\beta-1-k} u^{n-2}(x,t_0), \hspace{0.5em} \forall \beta \in \mathbb{Z}^+ -\{0\}, \hspace{0.5em} n>1 
\end{eqnarray}
Stating this problem in this way means that $g =u$ is the particular case when $n=1$ and  $i=0$. Some interests lay in the solutions of this equation for various values of $i$ and $n$ and understand the behavior of these solutions. Finally, there are no conditions on the boundaries in  the statement of \eqref{ProblemsetupB} as we are interested in the form of the general solutions such as decaying waves \cite{Petit}. 
\subsection{Beyond the families of Energy operators}\label{GenOperator2}
Looking at \eqref{dtppfwithApsdg2},  one can wonder what is the set of solutions $\mathcal{S}_2(\mathbb{R}^2 )$ for \eqref{ProblemsetupB}. If we define the subsets of $\mathbf{S}^{-}(\mathbb{R}^2)$ for $i$ in $\mathbb{Z}^+$:
\begin{eqnarray}
\mathbf{P}^i = & &\{  g \in \mathbf{S}^{-}(\mathbb{R}^2), \forall t \in \mathbb{R}^+, x_0  \in \mathbb{R}|  g(x_0,t) = \partial_t^i (\Psi_1^{t,+}(u)(x_0,t))^n,  \hspace{0.5em} u  \in \mathbf{S}^{-}(\mathbb{R}^2), \nonumber \\
 &&   \hspace{0.5em}n \in \mathbb{Z}^+ -\{0\} \}
\end{eqnarray} 
and,
\begin{eqnarray}
\mathbf{Q}^i = & & \{  g \in \mathbf{S}^{-}(\mathbb{R}^2), \forall x \in \mathbb{R}, t_0 \in \mathbb{R}^+| g(x,t_0) = \partial_x^i (\Psi_1^{t,+}(u)(x,t_0))^n,  \hspace{0.5em} u  \in \mathbf{S}^{-}(\mathbb{R}^2), \nonumber \\
&&  \hspace{0.5em} n \in \mathbb{Z}^+ -\{0\} \}
\end{eqnarray} 
with $\mathcal{S}_2(\mathbb{R}^2 )= \big ( \bigcup_{i \in \mathbb{Z}^+}  \mathbf{P}^i \big ) \bigcup \big ( \bigcup_{i \in \mathbb{Z}^+}  \mathbf{Q}^i \big )$. Note that one can define $\tilde{\mathcal{S}_2}(\mathbb{R}^2 )$ a similar space  as $\mathcal{S}_2(\mathbb{R}^2 )$ using the operator $\Psi_1^{x,+}$ instead of $\Psi_1^{t,+}$.
\newline We can decompose the successive derivatives of $(\Psi_1^{t,+}(u))^n$ (respectively $(\Psi_1^{x,+}(u))^n$) with the generalized energy operators, defined in \eqref{generalizedEO}, using Theorem $1$ with $u$ in $\mathbf{s}_1^{-}(\mathbb{R}^2)$  such as ( $x_0 \in \mathbb{R}$ , and $t_0 \in \mathbb{R}^+$):
\begin{eqnarray}\label{dtppfwithApsdg2}
\partial_t^{\beta} (\Psi_1^{t,+}(u)(x_0,t))^n &=& \sum_{k=0}^{\beta-1} \big(_{k}^{\beta-1} \big) \frac{n}{2}(\partial_t^k [[u(x_0,t)]^1]_1^{t,+} +\partial_t^k[[u(x_0,t)]^1]_1^{t,-} \partial_t^{\beta-1-k} u^{n-2}(x_0,t), \nonumber \\
& & \hspace{0.5em} \forall \beta \in \mathbb{Z}^+ -\{0\}, \hspace{0.5em} n>1 \\
\partial_x^{\beta} (\Psi_1^{t,+}(u)(x,t_0))^n &=& \sum_{k=0}^{\beta-1} \big(_{k}^{\beta-1} \big) \frac{n}{2}(\partial_x^k [[u(x,t_0)]^1]_1^{t,+} +\partial_x^k[[u(x,t_0)]^1]_1^{t,-} \partial_x^{\beta-1-k} u^{n-2}(x,t_0), \nonumber \\
& & \hspace{0.5em} \forall \beta \in \mathbb{Z}^+ -\{0\}, \hspace{0.5em} n>1
\end{eqnarray}
In the same way, we can define $\mathcal{S}_m(\mathbb{R}^2)$ $\subsetneq \mathbf{S}^{-}(\mathbb{R}^2)$ ($m$ in $\mathbb{Z}^+-\{0,1\}$)  associated with the generalized energy operators  $\big ([[.]^{m-2}]_k^{t,+}\big)_{k\in\mathbb{Z}}$ and $\big ([[.]^{m-2}]_k^{t,-}\big)_{k\in\mathbb{Z}}$,  and  $\tilde{\mathcal{S}_m}(\mathbb{R}^2)$ $\subsetneq \mathbf{S}^{-}(\mathbb{R}^2)$ ($m$ in $\mathbb{Z}^+-\{0,1\}$)  associated with the generalized energy operators  $\big ([[.]^{m-2}]_k^{x,+}\big)_{k\in\mathbb{Z}}$ and  $\big ([[.]^{m-2}]_k^{x,-}\big)_{k\in\mathbb{Z}}$. However, this model may define new sets of solutions for the linear PDEs thanks to the energy operators and the generalized energy operators. 
\section{Application to the homogeneous Helmholtz equation}
In the previous sections, we showed that it is possible to define the solutions of linear PDEs (e.g, \eqref{ProblemsetupB}) with the families of energy operator and to some extent the families of generalized energy operators. The theory is now applied to the particular case of the (homogeneous) Helmholtz equation.
%
\subsection{The homogeneous Helmholtz equation with solutions in the Schwartz space}
From \cite{JPMontillet2013} and \cite{Pain2005}, the homogeneous equation can be formulated  with $\alpha =2$, $h(x,t) =0$, $a_1 =1$ and $a_2= \frac{-1}{c^2}$ in \eqref{ProblemsetupB} such as:
\begin{equation}\label{eqpdep}
 \left\{
\begin{array}{rl}
\partial_x^{2}g(x,t) &  - \frac{1}{c^2} \partial_t^{2}g(x,t) = 0 ,\\
x\in \mathbb{R}, \hspace{0.5em}  x_0 \in \mathbb{R}, &  t\in \mathbb{R}^+,   \hspace{0.5em} t_0 \in \mathbb{R}^+\\
\end{array} \right.
\end{equation}
%
%
$c$ is the speed of light. It is well-known that the general solution $g(x,t)$ of this equation is a sum of two waves travelling in opposite direction such as $g(x,t) = u_1(t-x/c)+u_1(t+x/c)$ (e.g., \cite{Pain2005}). In this model,  $u_1(t-x/c)$ and $u_2(t-x/c)$ are the particular case of  $g(x,t) = \partial_t^i u^n(x,t)$ when $n=1$ and  $i=0$. A possible application of this theory is to look at the solutions for various values of $n$ and $i$.
 Applying the same development as in the previous section and looking for the solutions $g$ in $\mathcal{S}_1(\mathbb{R}^2)$, one can write \eqref{eqpdep} with the energy operator using Theorem $1$ (with $u$ in $\mathbf{s}^-(\mathbb{R}^2)$):
\begin{eqnarray}\label{eqpde2v}
\partial_x^{2}u^n(x,t) - \frac{1}{c^2} \partial_t^{2}u^n(x,t) &=&0 \nonumber \\
 - \frac{1}{c^2}\frac{n}{2}(\partial_t\Psi_1^{+,t}(u(x,t)) u^{n-2}(x,t) + \Psi_1^{-,t}(u(x,t)) \partial_t u^{n-2}(x,t))& &\nonumber\\
  + \frac{n}{2}(\partial_x\Psi_1^{+,x}(u(x,t)) u^{n-2}(x,t) + \Psi_1^{-,x}(u(x,t)) \partial_x u^{n-2}(x,t)) &=&0\nonumber\\
%
\end{eqnarray}
if $n=2$, \eqref{eqpdep}  can be simplified such as: 
\begin{eqnarray}\label{eqpde2v}
\partial_x^{2}u^2(x,t) - \frac{1}{c^2} \partial_t^{2}u^2(x,t) &=&0 \nonumber \\
- \frac{1}{c^2}(\partial_t\Psi_1^{+,t}(u(x,t)))+ (\partial_x\Psi_1^{+,x}(u(x,t)) &=& 0\nonumber\\
%
\end{eqnarray}
With the derivation chain rules property of the DEOs (e.g, \cite{JPMontillet2013}), we have  the equality $\frac{\partial (\Psi_k^{+,x}(g))}{\partial x} = \Psi_{k+1}^{+,x}(g) + \Psi_{k-1}^{+,x}(\partial_x g)$. Then, the previous equation becomes:
\begin{eqnarray}\label{eqpast}
\Psi_2^{+,x}(g)- \frac{1}{c^2} \Psi_2^{+,t}(g)  &=&0 \nonumber \\
\end{eqnarray}
\eqref{eqpast} agrees with the results previously published in \cite{JPMontillet2010}. Note that from \eqref{eqpde2v} when choosing some particular solutions $u$ of \eqref{eqpdep} (see next section), we can write $\Psi_1^{+,t}(u) \varpropto H \Psi_1^{+,x}(u)$ ($H$ in $\mathbb{R}$). Then, $\mathcal{S}_2(\mathbb{R}^2)$  and $\tilde{\mathcal{S}_2}(\mathbb{R}^2)$ are dual space.
\newline \indent Section \ref{sectionAppPDEs} recalls how the energy operators $(\Psi^{t,+}_k)_{k \in \mathbb{Z}}$ and $(\Psi^{t,-}_k)_{k \in \mathbb{Z}}$ decomposes uniquely $\partial_x^i u^n$ and $\partial_t^i u^n$. By estimating the energy operators  $(\Psi^{t,+}_k)_{k \in \mathbb{Z}}$, $(\Psi^{t,-}_k)_{k \in \mathbb{Z}}$, $(\Psi^{x,+}_k)_{k \in \mathbb{Z}}$ and $(\Psi^{x,-}_k)_{k \in \mathbb{Z}}$ for a given solution $u$ (in $\mathbf{s}^-(\mathbb{R})$), it is possible to estimate $k_1$ and $k_2$ in  $\mathbb{Z}^+$ such as for for $i > k_1$ $\mathbf{X}^i \sim \{0\}$ and for $j > k_2$ $\mathbf{Y}^j \sim \{0\}$. In other words, the energy operators should help to define a subset of solutions $\mathcal{A}(\mathbb{R}^2 ) \subseteq \mathcal{S}_1(\mathbb{R}^2 )$, such as it exists $k_1$ and $k_2$ in  $\mathbb{Z}^+$ for $\mathcal{A}(\mathbb{R}^2 )= \bigcup_{i \in [0,k_1]}  \mathbf{X}^i \bigcup \bigcup_{j \in [0,k_2]}  \mathbf{Y}^i$
\vspace{0.5em}
\newline Furthermore, one can use the generalized energy operator families to define the solutions of the Helmholtz equation in $\mathcal{S}_2(\mathbb{R}^2)$ defined in Section \ref{GenOperator2}. Thus, one can write for $g$ in $\mathcal{S}_2(\mathbb{R}^2)$ (and $u$ in $\mathbf{s}_1^{-}(\mathbb{R}^2) \subsetneq \mathbf{S}^{-}(\mathbb{R}^2)$):
%
\begin{eqnarray}\label{eqpde2vv}
\partial_x^2 (\Psi_1^{+,t}(u))^2 - \frac{1}{c^2} \partial_t^2 (\Psi_1^{+,t}(u))^2  &=&0 \nonumber \\
\partial_x^2([[u(x,t)]^0]_1^{t,+})^2 - \frac{1}{c^2} \partial_t^2([[u(x,t)]^0]_1^{t,+})^2 &=&0 \nonumber \\
\partial_x\big ([[u(x,t)]^1]_1^{t,+}+[[u(x,t)]^1]_1^{t,-} \big )- \frac{1}{c^2} \partial_t^2\big ([[u(x,t)]^1]_1^{t,+}+[[u(x,t)]^1]_1^{t,-} \big ) &=&0 \nonumber \\
\end{eqnarray}
Furthermore, we can define other sets of solutions in $\mathbf{S}^-(\mathbb{R}^2)$ similar to the definitions of $\mathcal{S}_1(\mathbb{R}^2)$ in Section \ref{sectionAppPDEs} and  $\mathcal{S}_2(\mathbb{R}^2)$ in Section  \ref{GenOperator2} with other families of generalized energy operators such as $\big (([[.]^2]_k^{x,+}\big)_{k\in\mathbb{Z}}$, $\big (([[.]^2]_k^{t,+}\big)_{k\in\mathbb{Z}}$.
\vspace{0.5em}
\newline \indent Finally, can we consider  $\big (([[u]^m]_k^{x,+}\big)_{k\in\mathbb{Z}}$, $\big (([[u]^m]_k^{t,+}\big)_{k\in\mathbb{Z}}$, $\big (([[u]^m]_k^{x,-}\big)_{k\in\mathbb{Z}}$, $\big (([[u]^m]_k^{t,-}\big)_{k\in\mathbb{Z}}$ ($m$ in $\mathbb{Z}^+$) as "wave" which are also propagating in time and space? They are generally not the full solutions of the homogeneous Helmholtz equation (see Section \ref{sectionAppPDEs}). Let us look at a case study. 
%
%
%
%
\subsection{Application to the case of the  evanescent wave}
According to the definition of the solutions of the Helmholtz equation given in the previous section for the case when the PDE is solved  in $S^-(\mathbb{R}^2)$, we choose here the solutions of the type:
\begin{equation}\label{evanescent}
 \left\{
\begin{array}{rl}
u(x,t) = & Real \{A \exp{(-k_1 x) } \exp{(j(\omega t -k_2 x))} \} ,\\
t \in [0, T ], \hspace{0.5em} x \in [x_1, x_2], & \hspace{0.5em} (x_1,x_2) \in \mathbb{R}^2, \hspace{0.5em} x_1<x_2\\
\end{array} \right.
\end{equation}
$k_1$ and $k_2$ are the wave numbers,  $\omega$ is the angular frequency and $A$ is the amplitude of this wave \cite{Petit}. Assuming $\omega$ and $k_2$ known, one can add some conditions to the limits in order to estimate $k_1$ and $A$. However, our interest is just the general form assuming that all the parameters are known.
\newline In addition, a travelling wave solution of \eqref{eqpdep} should satisfy some relationship between $k_1$, $k_2$ and $\omega$ \cite{Petit}. If we replace the general solution $g$ of \eqref{eqpdep} with solutions in $\mathcal{S}_1(\mathbb{R}^2)$, $\mathcal{S}_2(\mathbb{R}^2)$ or $\mathcal{S}_m(\mathbb{R}^2)$ ($m$ in $\mathbb{Z}^+-\{0,1\})$, one can show that we have always the relationship $(k_1+ j k_2)^2-(\frac{j\omega}{c})^2=0$ (or $Real\{(k_1+ j k_2)^2-(\frac{j\omega}{c})^2\}=0$).
\newline \indent It is also possible to estimate the average power of the wave through a section (S) such as (\cite{Amzallag et al.}, p. 74):
\begin{equation}
P_a = \frac{w}{T} = \frac{\int_S|g|^2dS}{T} 
\end{equation}
replacing with the evanescent waves (for $(a,b,L,t_0)$ in $\mathbb{R}^4$, $k$ in $\mathbb{Z}^+$, $n$ in $\mathbb{Z}^+-\{0\}$), 
\begin{eqnarray}
P_a &=& \frac{\int_a^{b}|\partial_t^{k} u^n(x,t_0)|^2 Ldx}{T} \nonumber \\
P_a &=& \frac{\int_a^{b}Real\{ (nj\omega)^{2(k+1)}u^{2n}(x,t_0)\} Ldx}{T} \nonumber \\
P_a &=& \frac{1}{T}  Real\{ \frac{1}{2n(-k_1-jk_2)}L(nj\omega)^{2(k+1)}[u^{2n}(x,t_0)]_a^b\}
\end{eqnarray}
Thus, comparing $P = \frac{\int_a^{b}|u(x,t_0)|^2 Ldx}{T}$ with $P_a$, we can write $P_a = \alpha P$ ($\alpha$ in $\mathbb{R}$). We can also estimate the power $P_a$ with $g$  in $\mathcal{S}_2(\mathbb{R}^2)$ or $\mathcal{S}_m(\mathbb{R}^2)$. This will change the value of $\alpha$. Thus, it turns out that we are "demultiplying" the averaged power.
\newline \indent Finally, this work ends with the computation of \eqref{evanescent} in $2D$ for some solutions in $\mathcal{S}_1(\mathbb{R}^3)$ and $\mathcal{S}_2(\mathbb{R}^3)$ in $\mathbf{S}^-(\mathbb{R}^3)$. The figures below show some examples with $\partial_t^3u^n(x,t)$ and $\partial_t^3\big ( \Psi_1^{t,+} (u (x,t)) \big)^n$ ($n$ in $\{2,5,8\}$). Note that $k_1$ is equal to $-50$ $cm^{-1}$, $k_2$ is equal to $100$ $cm^{-1}$ and $A$ equal to $10$ $cm$. The values are chosen arbitrarily. One can see that for fixed values of $t$, $x$ and $y$,  the wave is contracting when using higher and higher values of $n$. This result is expected, because when $n$ increases the attenuation coefficient in \eqref{evanescent} is equal to $\exp{(-n k_1 x) }$. 
\section{Conclusions}
This work generalizes the $\bold{Lemma}$ $0$ and $\bold{Theorem}$ $0$ shown in \cite{JPMontillet2013} using the families of generalized energy operators $\big ([[.]^p]_k^+\big)_{k\in\mathbb{Z}}$ and $\big ([[.]^p]_k^- \big)_{k\in\mathbb{Z}}$ ($p$ in $\mathbb{Z}^+$) defined in Section \ref{sectionenergyopdefandgen} (e.g, \eqref{bracketl} and \eqref{eqpsip2}). $\bold{Lemma}$ $1$ shows that the successive derivatives of $\big ([[f]^{p-1}]_1^+ \big)^n$ ($n$ in $\mathbb{Z}^+$, $n>1$) can be decomposed with the generalized energy operators  $\big ([[.]^p]_k^+\big)_{k\in\mathbb{Z}}$ when $f$ is in the subspace $\mathbf{S}_p^-(\mathbb{R})$. With $\bold{Theorem}$ $1$ and $f$ in $\mathbf{s}_p^-(\mathbb{R})$, one can decompose uniquely the successive derivatives of $\big ([[f]^{p-1}]_1^+ \big)^n$ ($n$ in $\mathbb{Z}^+$, $n>1$) with the generalized energy operators  $\big ([[.]^p]_k^+\big)_{k\in\mathbb{Z}}$ and $\big ([[.]^p]_k^-\big)_{k\in\mathbb{Z}}$. $\mathbf{S}_p^-(\mathbb{R})$ and $\mathbf{s}_p^-(\mathbb{R})$ ($p$ in $\mathbb{Z}^+$) are subspaces of the Schwartz space $\mathbf{S}^-(\mathbb{R})$. Their definitions involve the so called energy spaces defined in $\bold{Definition}$ $3$. 
 The proofs of $\bold{Lemma}$ $1$ and $\bold{Theorem}$ $1$ follow a similar structure : an induction on both $p$ and $n$. Note that the special case $ n <-1$ and $n= \pm 1$ are discussed at the end of Section \ref{sectionenergyopdefandgen}. It is worth emphasizing that in the case $p=0$, $\bold{Lemma}$ $1$ and $\bold{Theorem}$ $1$ are the same statements as $\bold{Lemma}$ $0$ and $\bold{Theorem}$ $0$. This demonstrates that this work generalizes the previous work of \cite{JPMontillet2013}.
\newline The second part of the work focuses on the application of the theory developed in \cite{JPMontillet2013} and with the generalized energy operators  to the linear PDEs. In this application, the solutions are functions of two variables defined onto $\mathcal{S}_m(\mathbb{R}^2)$ $\subsetneq \mathbf{S}^{-}(\mathbb{R}^2)$ ($m$ in $\mathbb{Z}^+-\{0,1\}$) $\mathbf{S}^-(\mathbb{R}^2)$. The definition of these subspaces include the energy operators and the generalized energy operators using $\mathbf{Theorem}$ $0$ and $\mathbf{Theorem}$ $1$. These subspaces define a  new model for the solutions of the linear PDEs on to $\mathbf{S}^-(\mathbb{R}^2)$. We show the application of this model in the case of the Helmholtz equation and the particular case of evanescent waves.  This opens possible applications of this theory in astronautic and astrophysics in "demultiplying" the averaged power of the energy operators and generalized energy operators (e.g. \cite{Labun}).
\section*{Acknowledgment}
A special thanks is addressed to Professor Alan McIntosh at the Centre for Mathematics and its Applications at the Australian National University (ANU) for its inputs and discussions when writing this manuscript. The author also acknowledges the comments from Dr. Igor Ivanov from the Atomic and Molecular Physics laboratory at the ANU, and Dr. Malcolm S. Woolfson from the School of Electrical Engineering at the University of Nottingham (UK).
%
\newpage
\begin{figure}
\begin{tabular}{cc}
\hspace{-8em}\subfloat[]{\includegraphics[width=3.4in]{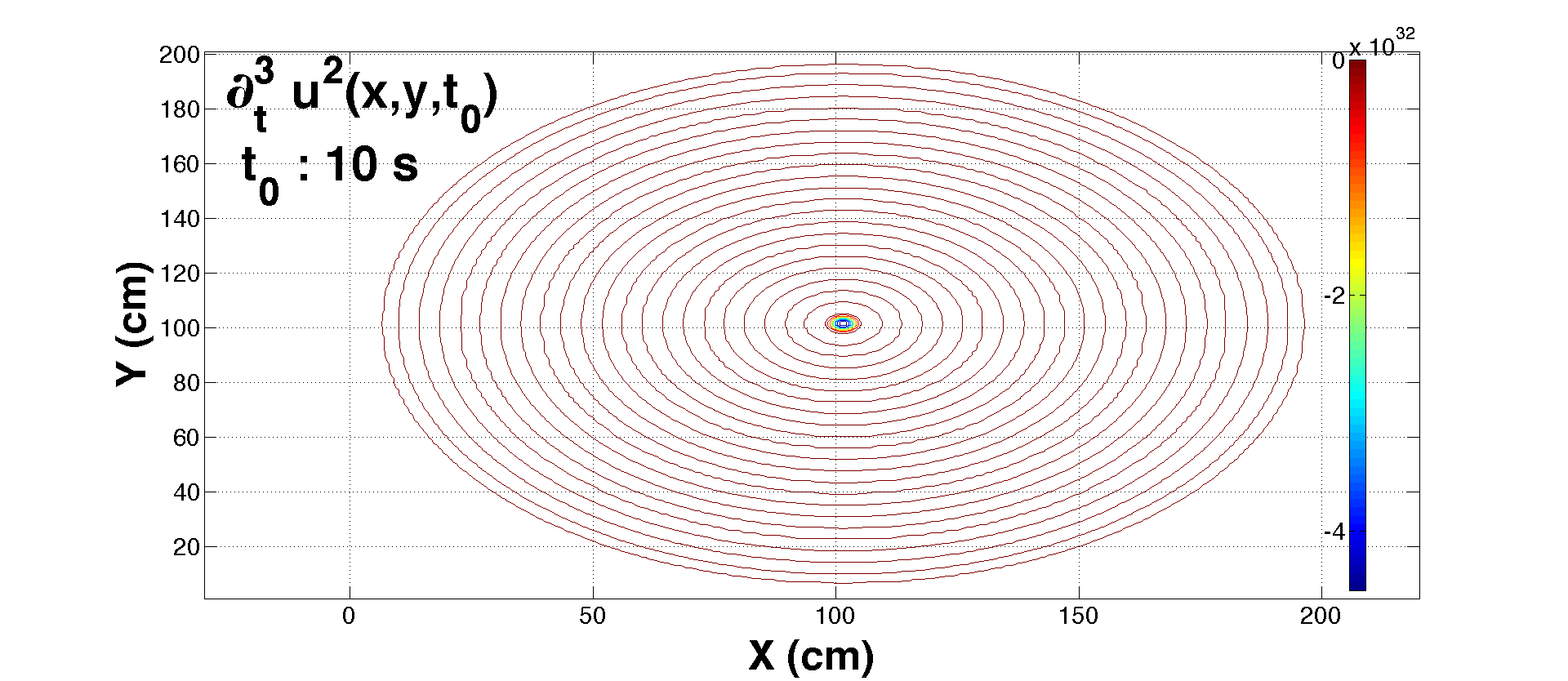}} &
\hspace{-1em}\subfloat[]{\includegraphics[width=3.4in]{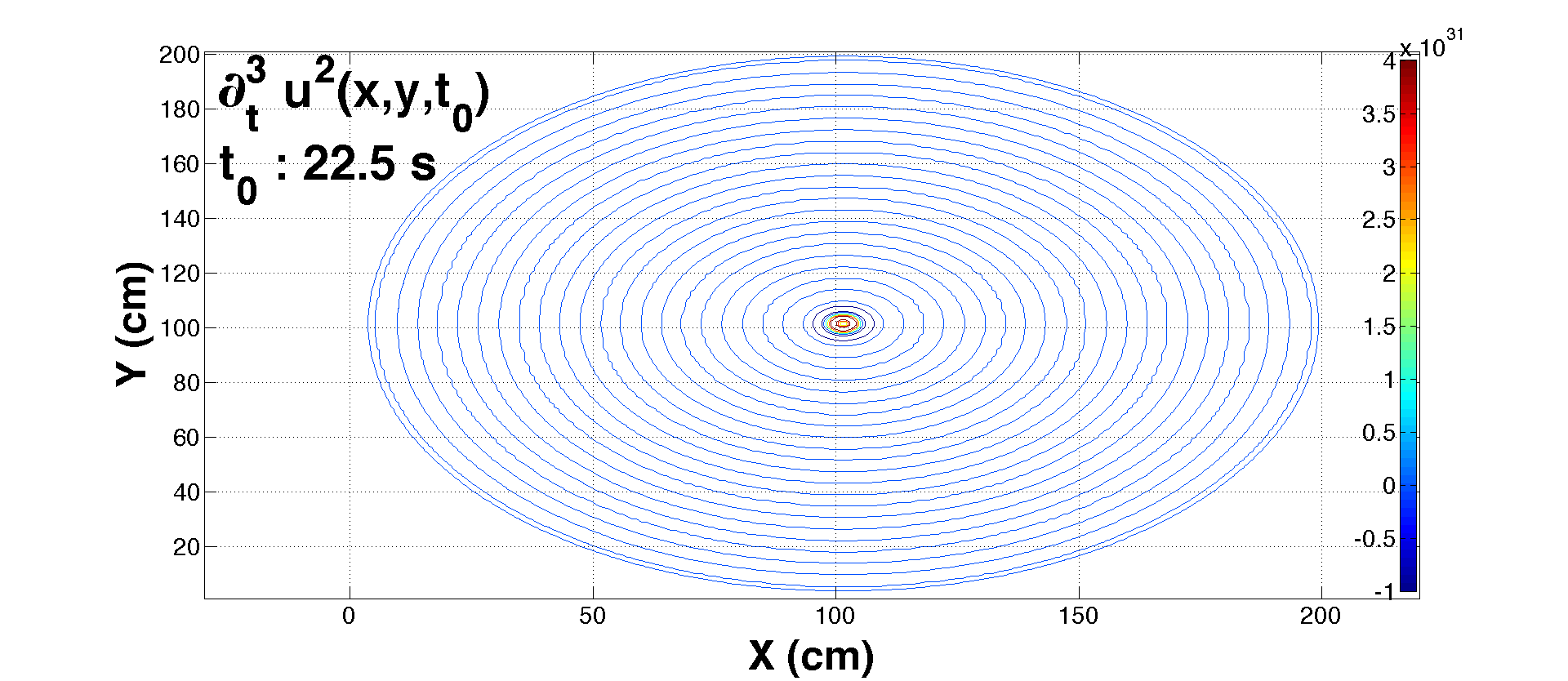}} \\
\hspace{-8em}\subfloat[]{\includegraphics[width=3.4in]{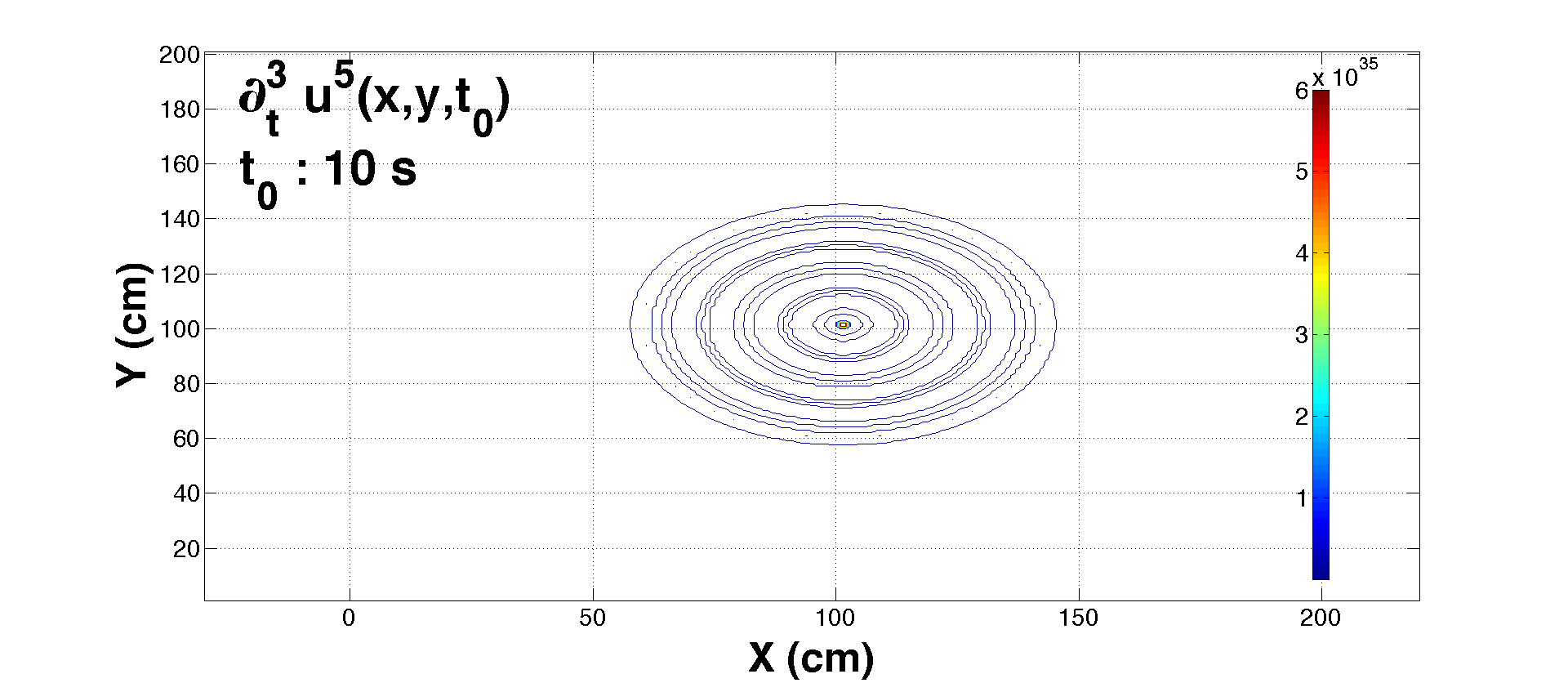}} &
\hspace{-1em}\subfloat[]{\includegraphics[width=3.4in]{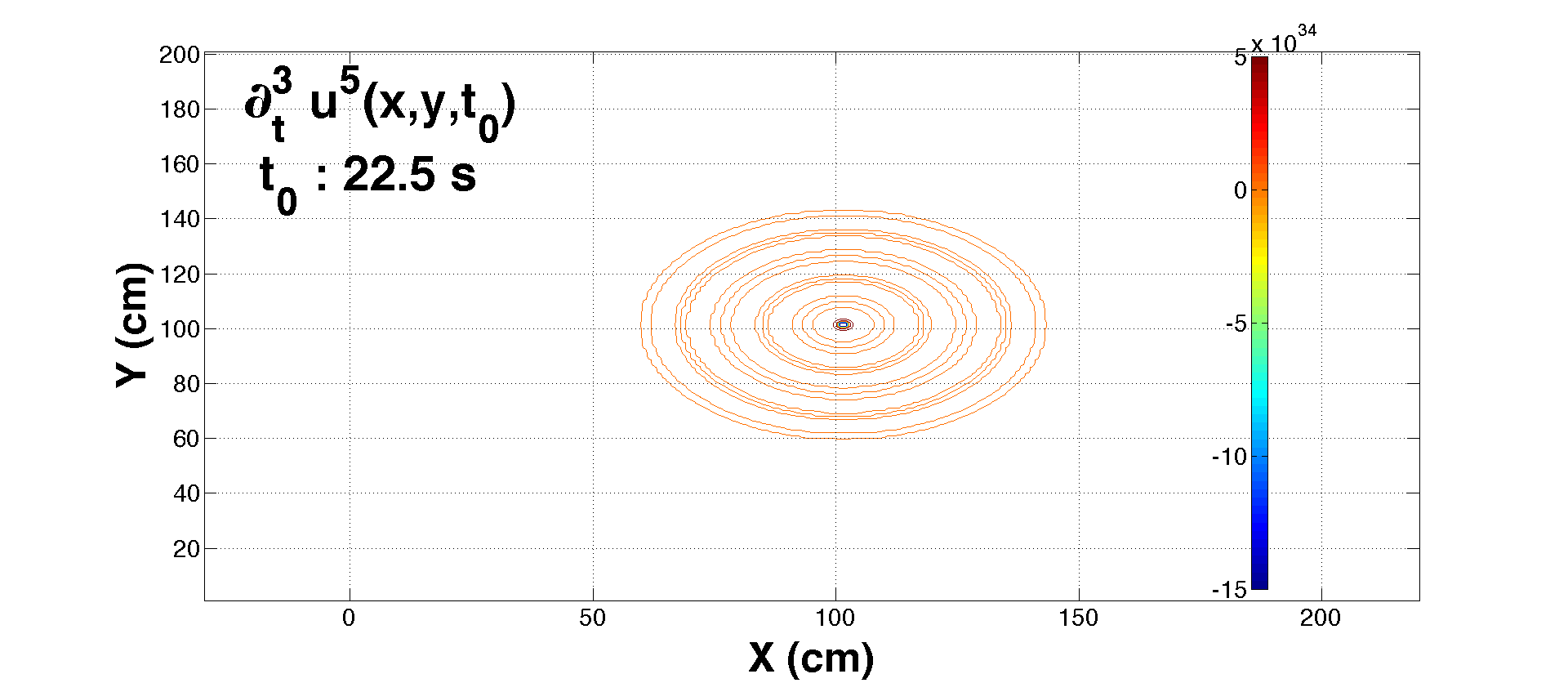}}  \\
\hspace{-8em}\subfloat[]{\includegraphics[width=3.4in]{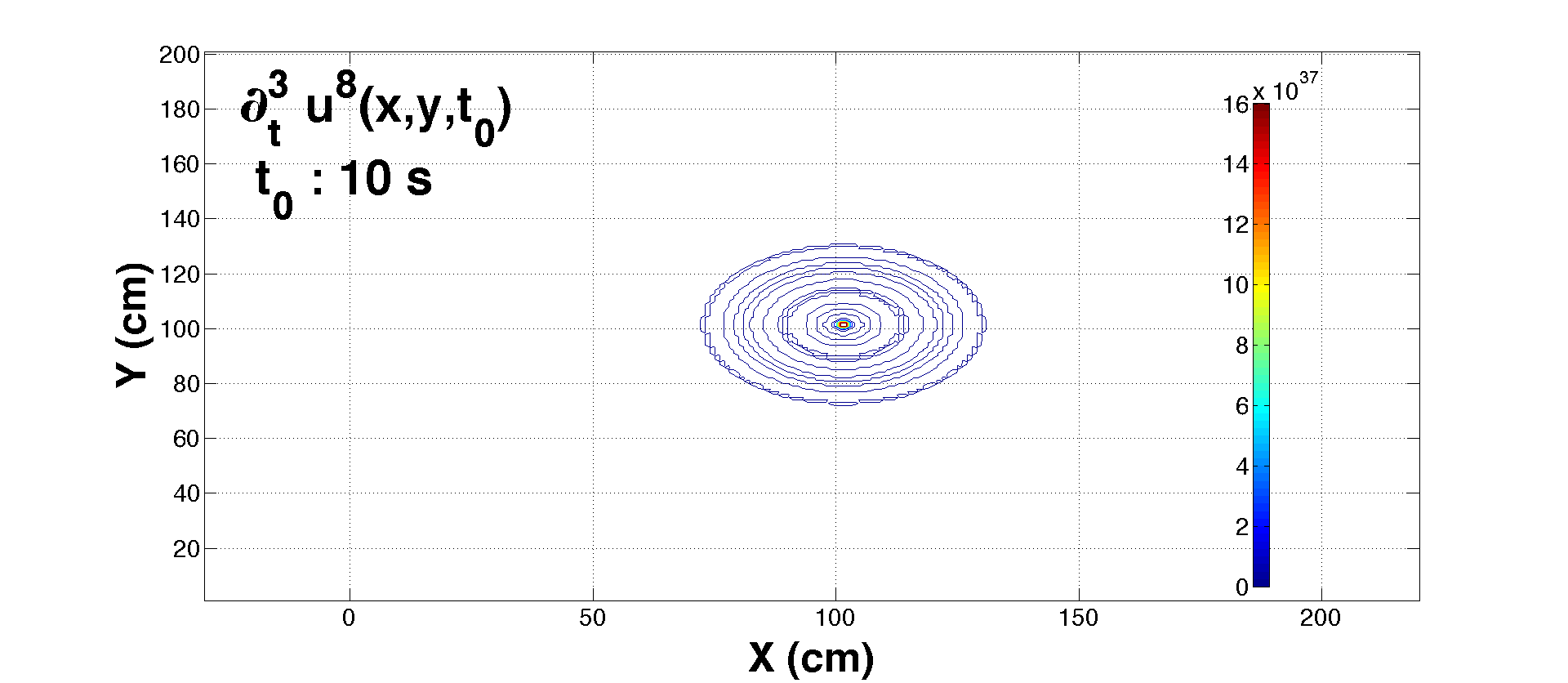}} &
\hspace{-1em}\subfloat[]{\includegraphics[width=3.4in]{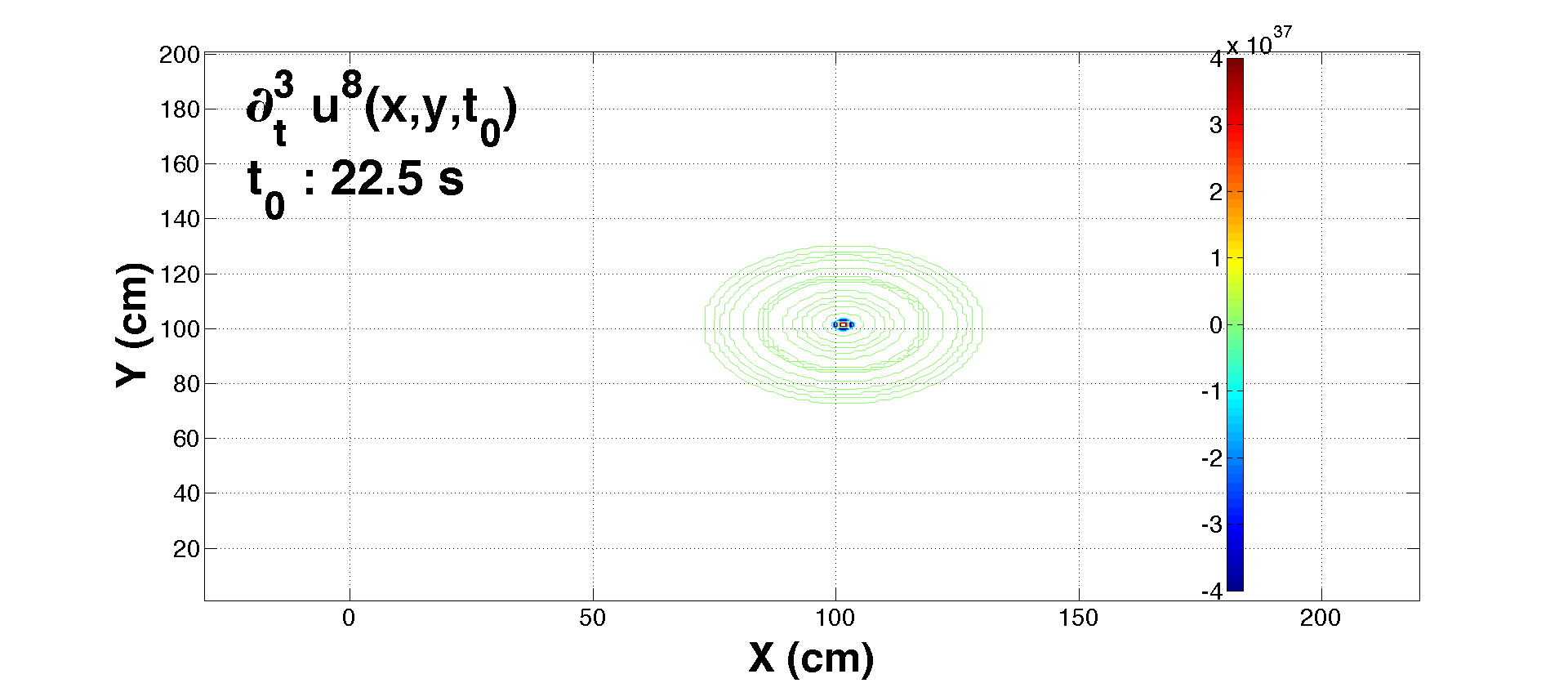}} 
\end{tabular}\label{Figure1}
\caption{Computation of $\partial_t^3u^n(x,t)$ ($n$ in $\{2,5,8\}$)}
\end{figure}
\begin{figure}\label{Figure2}
\begin{tabular}{cc}
\hspace{-8em}\subfloat[]{\includegraphics[width=3.4in]{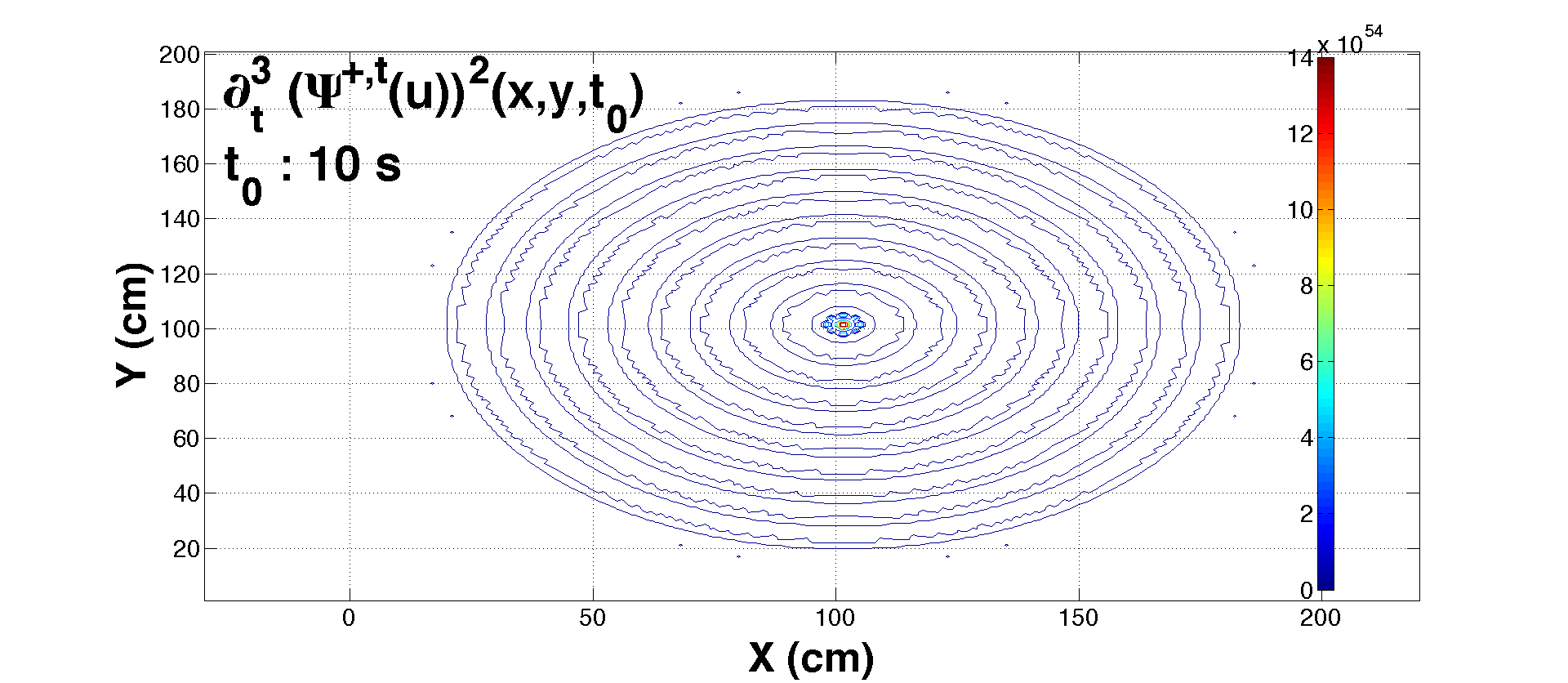}} &
\hspace{-1em}\subfloat[]{\includegraphics[width=3.4in]{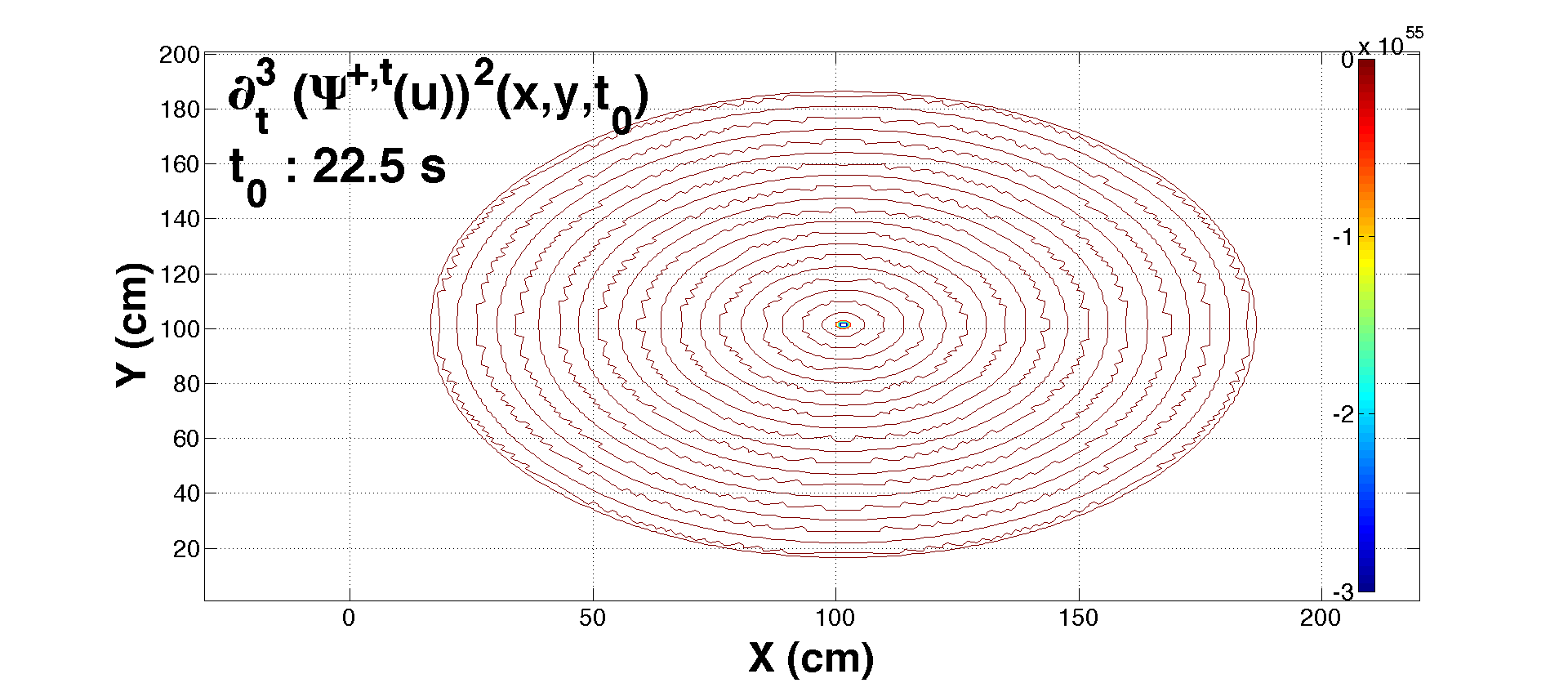}} \\
\hspace{-8em}\subfloat[]{\includegraphics[width=3.4in]{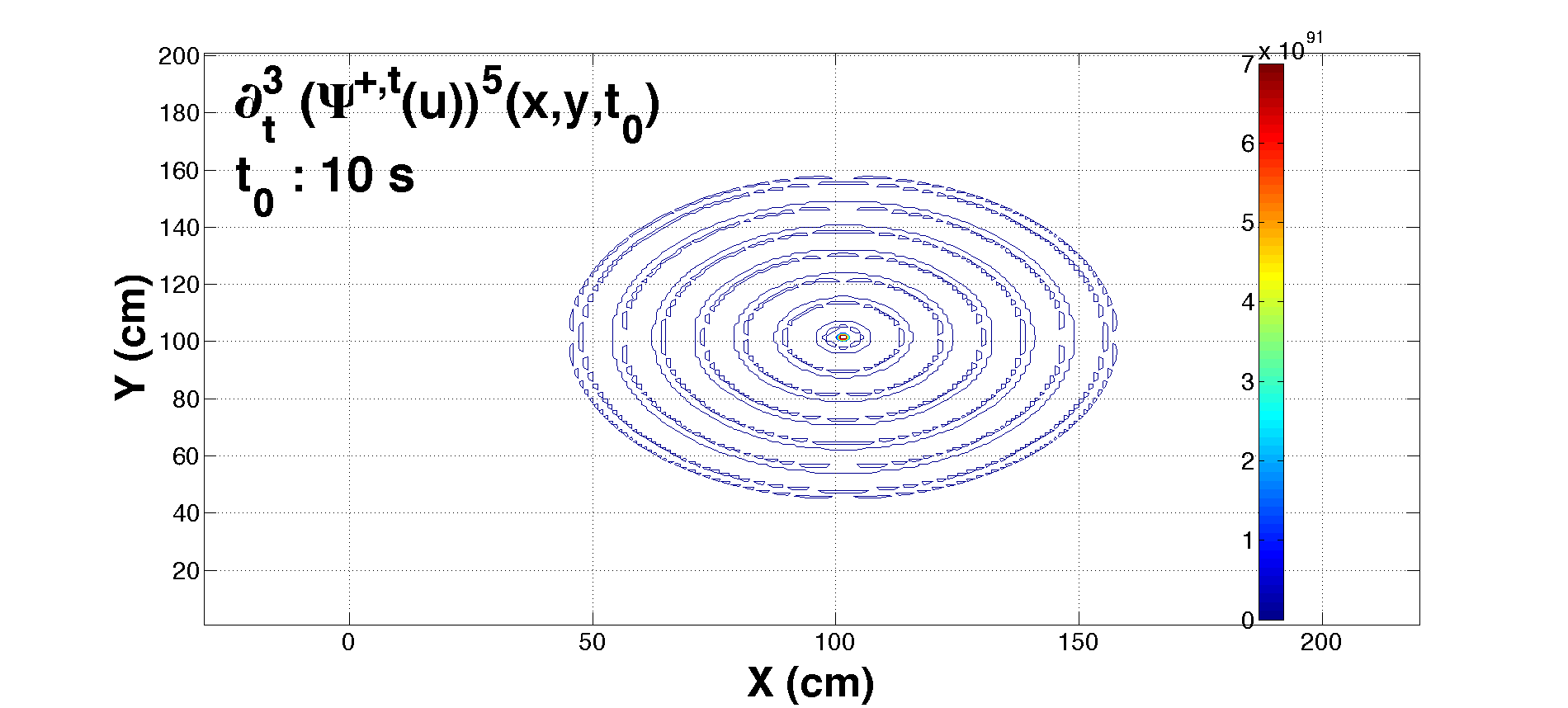}} &
\hspace{-1em}\subfloat[]{\includegraphics[width=3.4in]{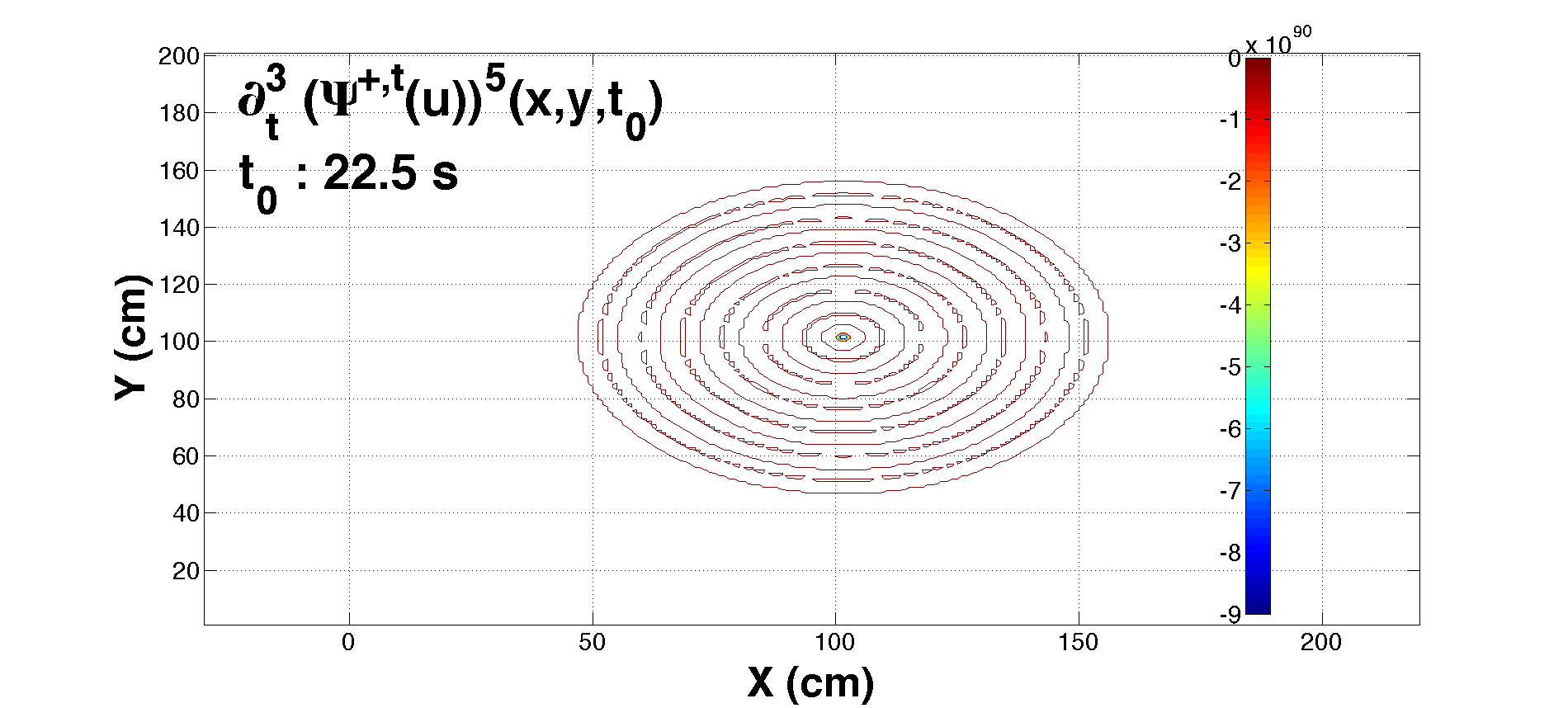}}  \\
\hspace{-8em}\subfloat[]{\includegraphics[width=3.4in]{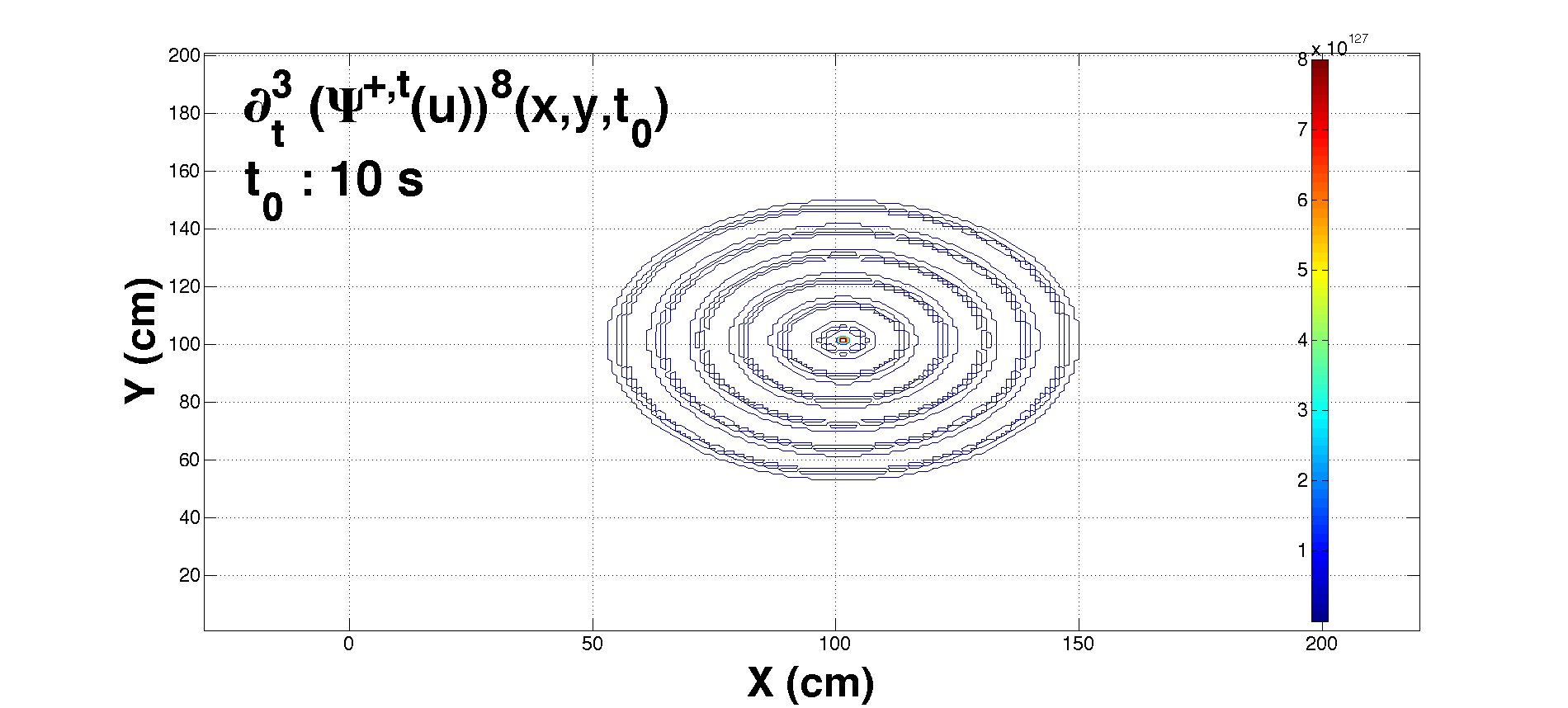}} &
\hspace{-1em}\subfloat[]{\includegraphics[width=3.4in]{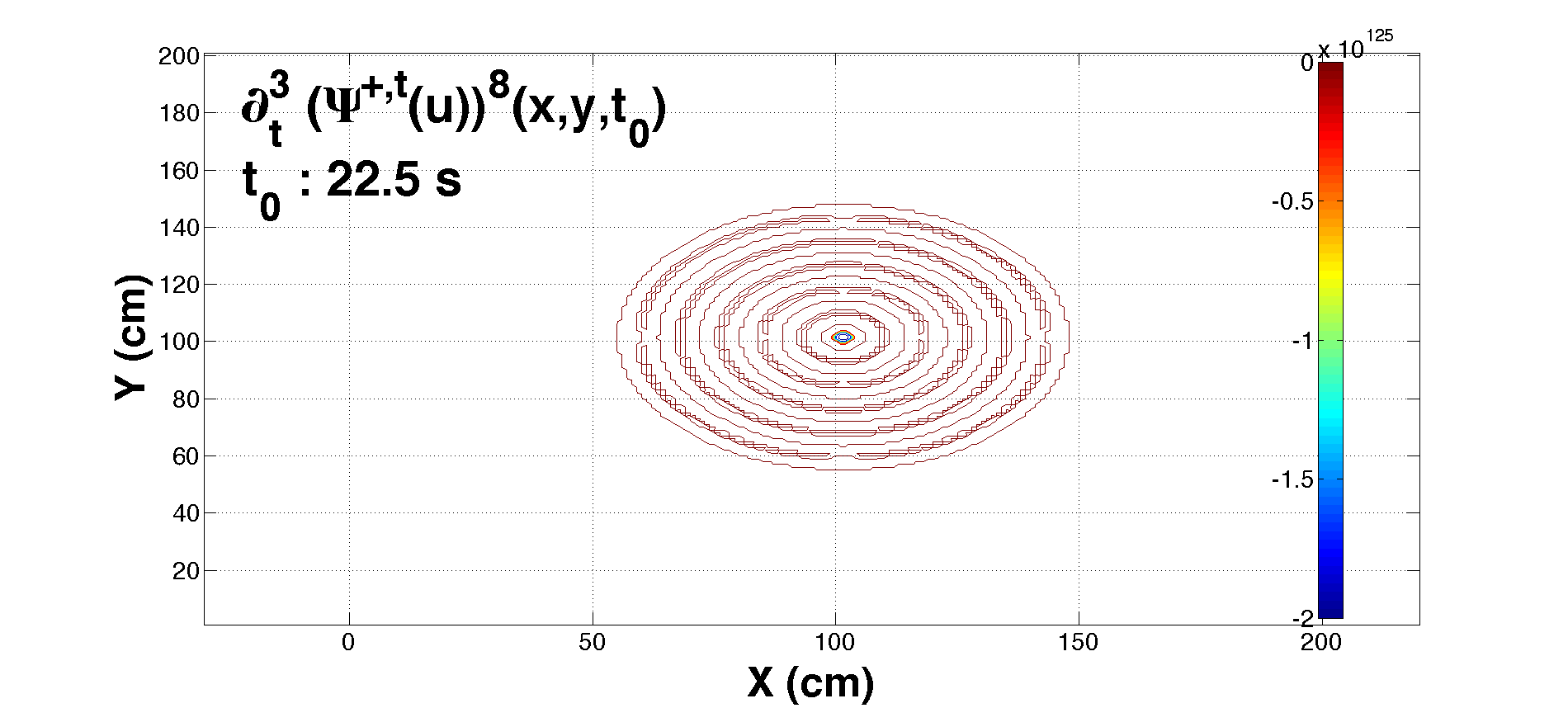}} 
\end{tabular}
\caption{Computation of $\partial_t^3\big ( \Psi_1^{t,+} (u (x,t)) \big)^n$ ($n$ in $\{2,5,8\}$)}
\end{figure}
\end{document}